
\input amstex
\documentstyle{amsppt}

\magnification=1200

\TagsOnRight
\loadbold

\redefine\.{\kern1pt.}
\define\mz{\<${}'$\<}
\redefine\<{\leavevmode\kern-\mathsurround}
\redefine\le{\leqslant}
\redefine\ge{\geqslant}
\define\({\hskip1pt(}
\define\>{\hskip1pt}
\define\){)\hskip1pt}
\define\mi{\kern0.8pt|\kern0.8pt}
\define\op{\operatorname}
\define\ov{\overline}
\define\bdot{\boldsymbol\cdot}

\NoRunningHeads  
\topmatter  
\title Ergodic unitarily invariant measures on the space\\   
of infinite Hermitian matrices\endtitle  
\endtopmatter  
  
\document  
  
\centerline{{\smc Grigori Olshanski}}  
\medskip {\eightpoint  
\centerline{Institute for Problems of Information Transmission}   
\centerline{Bolshoi Karetnyi per. 19, 101447 Moscow GSP--4, RUSSIA}  
\centerline{E-mail: olsh$\@$ippi.ac.msk.su}}  
\medskip  
\centerline{{\smc and}}  
\medskip  
\centerline{{\smc Anatoli Vershik}}  
\medskip {\eightpoint  
\centerline{St.Petersburg Branch of Steklov Mathematical Institute}  
\centerline{Fontanka 27,  191011 St.Petersburg, RUSSIA}   
\centerline{E-mail: vershik$\@$pdmi.ras.ru}}  
\bigskip    
\centerline{January 1996}
\bigskip
To appear in F.~A.~Berezin's memorial
volume, American Mathematical Society, 1996.

\newpage
\centerline{{\smc Abstract}}  
\bigskip  
  
Let $H$ be the space of all Hermitian matrices of infinite
order and $U(\infty)$ be the inductive limit of the chain $U(1)\subset
U(2)\subset\ldots$ of compact unitary groups. The group $U(\infty)$
operates on the space $H$ by conjugations, and our aim is to classify
the ergodic $U(\infty)$-invariant probability measures on $H$ by
making use of a general asymptotic approach proposed in Vershik's
note \cite{V}. The problem is reduced to studying the limit behavior
of orbital integrals of the form
$$
\int_{B\in\Omega_n}e^{i\op{tr}(AB)}M_n(dB),
$$
where $A$ is a fixed $\infty\times\infty$ Hermitian matrix with
finitely many nonzero entries, $\Omega_n$ is a $U(n)$-orbit in the
space of $n\times n$ Hermitian matrices, $M_n$ is the normalized
$U(n)$-invariant measure on the orbit $\Omega_n$, and $n\to\infty$.

We also present a detailed proof of an ergodic theorem for inductive
limits of compact groups that has been announced in \cite{V}.

There is a remarkable link between our subject and Schoenberg's
\cite{S2} theory of totally positive functions, and our approach
leads to a new proof of Schoenberg's \cite{S2} main theorem,
originally proved by function-theoretic methods.

On the other hand, our results have a representation-theoretic
interpretation, because the ergodic $U(\infty)$-invariant measures on
$H$ determine irreducible unitary spherical representations of an
infinite--dimensional Cartan motion group.

The present paper is closely connected with a series of articles by
S. V. Kerov and the authors on the asymptotic representation theory
of ``big'' groups, but it can be read independently. 
\footnote""{ Both authors were  
partially supported by the International Science Foundation under Grant  
MQV000. The first--named author was also partially supported by the
Russian Foundation for Basic Research under Grant 95--01--00814.}
    
\newpage  
\hfill {\it To the memory of F. A. Berezin}  
\bigskip
\head\S0.\enspace Introduction \endhead

\subhead\rom{1.}\quad Ergodic measures \endsubhead The description of ergodic
invariant measures for group actions is a traditional problem of
ergodic theory. It is well known that it is not always that this
problem can be solved in a satisfactory way. It is a surprising fact that for
certain infinite-dimensional (or ``big''\footnote{The term ``big
groups'' was suggested by one of us in \cite{V}. This term has no
rigorous definition, it can mean ``infinite-dimensional'' or ``out
of the class of locally compact groups'' or something like that. For
instance, the group $S(\infty):=\varinjlim S(n)$, the inductive limit
of finite symmetric groups, is a discrete, hence a locally compact
group, but its properties are very similar to that of the
infinite-dimensional group $U(\infty):=\varinjlim U(n)$, so we
prefer to rank $S(\infty)$ among the ``big'' groups.}) groups, there
exist nice actions whose ergodic measures can be completely
described. The simplest examples are given by the classical theorems
due to B.~de~Finetti and I.~J.~Schoenberg.

We recall that, by de Finetti's theorem, the ergodic measures on a
product space $X^{\infty}$, invariant under the group of permutations
of the copies of $X$ (the infinite symmetric group), are exactly the
product measures with identical factors. Schoenberg's theorem
\cite{S1} states that the ergodic $O(\infty)$-invariant measures on
the space ${\Bbb R}^{\infty}$ are exactly the Gaussian product measures
(by $O(\infty)$ we denote the inductive limit group $\varinjlim
O(n)$).\footnote{About this theorem, see also
Berg--Christensen--Ressel \cite {BCR}.}

These classical examples and a number of more complicated ones were
discussed in the note \cite{V} by one of the authors. In that note, a
general ``ergodic method'' for inductive limits of compact groups was
proposed. This method was further developed in the papers \cite{VK1,
VK2, KV} by Vershik and Kerov. The aim of the present
paper is to give a detailed exposition of the method on a model
example considered in \cite{V}.

Let $H(n)$ denote the space of $n\times n$ complex Hermitian
matrices, and let $H=\varprojlim H(n)$ be the space of all infinite
Hermitian matrices. Let $U(\infty)=\varinjlim U(n)$ be the group of
infinite unitary matrices $u=[u_{ij}]$ such that $u_{ij}=\delta_{ij}$
when $i+j$ is large enough. The group $U(\infty)$ operates on the
space $H$ by conjugations, and we are interested in the class ${\Cal
M}$ of all ergodic $U(\infty)$-invariant Borel probability measures
on $H$.

Further, let $H(\infty)=\varinjlim H(n)$ be the space of
$\infty\times\infty$ Hermitian matrices with finitely many nonzero
entries. The spaces $H(\infty)$ and $H$ are in a natural duality, and
any measure on $H$ is uniquely determined by its characteristic
function (Fourier transform), which is a function on $H(\infty)$.

\proclaim{Classification Theorem} The characteristic functions of the
measures $M\in{\Cal M}$ are exactly those of the form
$$
f(A)=e^{i\gamma_1\op{tr}A-\gamma_2\op{tr}(A^2)/2}
\det\bigg(\prod_{k=1}^{\infty}\frac{e^{-ix_kA}}{1-ix_kA}\bigg),\qquad
A\in H(\infty),\tag0.1
$$
where $\gamma_1,\gamma_2,x_1,x_2,\dots$, the parameters of the
measure $M$, are real numbers such that $\gamma_2\ge0$ and $\sum
x_k^2<\infty$.
\endproclaim

The Classification Theorem implies that any measure $M\in{\Cal M}$
can be written as the convolution product of a Gaussian measure and a
(finite or countable) family of non-Gaussian ``elementary'' ergodic
measures. The latter are essentially supported by rank one
matrices and are related to Wishart distributions, well-known in
multivariate statistical analysis. Let us emphasize that in this
theorem, the structure of the answer turns out to be much more
complicated than in de Finetti's and Schoenberg's \cite{S1} theorems:
up to trivial exceptions, the ergodic measures $M\in{\Cal M}$ are
neither product nor Gaussian measures.

Note that formula (0.1) can be rewritten in the following form:
$$
f(A)=\prod_{a\in\op{Spec}(A)}F(a),\qquad A\in H(\infty),\tag0.2
$$
where $\op{Spec}\(A)$ is the collection of eigenvalues of $A$
(taken with their multiplicities) and
$$
F(a)=e^{i\gamma_1a-\gamma_2a^2\!/2}
\prod_{k=1}^{\infty}\frac{e^{-ix_ka}}{1-ix_ka},\qquad a\in{\Bbb R}\.\tag0.3
$$
The function $F(a)$ has a simple meaning. Let us regard the matrices
$B\in H$ as random matrix variables defined on the probability space
$(H,M)$. Then $F(a)$ is the characteristic function for the
distribution of any diagonal entry $B_{ii}$,
$i=1,2,\dots$.\footnote{This implies that the whole distribution of
the matrix entries is uniquely determined by the distribution of any
diagonal entry! Note that this property holds for ergodic measures $M$
only: it follows from the so-called multiplicativity theorem, see
Theorem 2.1 below.}

As shown by Pickrell \cite{Pi2}, the Classification Theorem
can be derived from a deep function-theoretic result due to
Schoenberg \cite{S2}. In the present paper, we prove the
Classification Theorem in an entirely different way: our method
consists in studying how an ergodic measure $M\in{\Cal M}$ is
approximated by finite-dimensional ``orbital measures''. (By an
{\it orbital measure\/} we mean the $U(n)$-invariant probability measure
supported by a $U(n)$-orbit $\Omega_n\subset H(n)$, where
$n=1,2,\dots$.)

We would like to emphasize that such an approach not only leads to
the description (0.1) of ergodic measures, but also provides us with
additional information about them. Namely, we find a
characterization of those sequences $\{M_n\}$ of orbital measures
that weakly converge, as $n\to\infty$, to an ergodic measure
$M\in{\Cal M}$, so we can understand how $M$ ``grows up'' from
orbital measures. In particular, we can relate the parameters
$\gamma_1,\gamma_2,x_1,x_2,\dots$ to the asymptotic behavior
of the eigenvalues of large Hermitian matrices.

A characterization of the convergent sequences $M_n$ is the main
result of the paper---Theorem 4.1. The statement is too long to be
reproduced here in detail, but roughly speaking, the picture is as
follows. Given $n=1,2,\dots$, pick any matrix from the orbit
$\Omega_n$ supporting $M_n$ and represent the collection of the
eigenvalues of that matrix as an $n$-point configuration on the real
line. Then this configuration, being contracted with scaling factor
$1/n$, must converge, as $n\to\infty$, to a (countable) point
configuration $\{x_1,x_2,\dots\}\subset{\Bbb R}$. (Here the $x_k$'s
coincide with the parameters appearing in (0.1); in fact, there are
also additional restrictions on the growth of the orbits related to
the parameters $\gamma_1$ and $\gamma_2$.)

\subhead\rom{2.}\quad Relationships with representations and totally positive
functions\endsubhead We are interested in the measures of
class ${\Cal M}$ for several reasons:

1) (The initial motivation of \cite{V}.) The action of the group
$U(\infty)$ on the space $H$ is a natural example of a ``big'' group
action. As compared with the action of the group $O(\infty)$ on the
space ${\Bbb R}^{\infty}$ (Schoenberg's case \cite{S1}) or the action
of $U(\infty)$ on ${\Bbb C}^{\infty}$, this example is of the next
level of complexity. So, from the viewpoint of ergodic theory, it is
interesting to find ergodic measures and to compare the result with
the one in Schoenberg's case.

2) The ergodic measures of class ${\Cal M}$ determine some
irreducible unitary representations of a certain
infinite-dimensional group $G(\infty)$.

3) These measures are closely related to Schoenberg's totally
positive functions.

4) They are also related to harmonic analysis on the group
$U(\infty)$ (cf.~\cite{KOV}).

We shall briefly comment on items 2) and 3).

The group $G(\infty)$ is a model example of an infinite-dimensional
Cartan motion group. It can be defined as the inductive limit
$G(\infty)=\varinjlim G(n)$, where $G(n)$ stands for the semidirect
product $U(n)\ltimes H(n)$. There is a natural bijective
correspondence between ergodic $U(\infty)$-invariant probability
measures on $H$ and irreducible unitary representations of the group
$G(\infty)$ that are spherical with respect to the subgroup
$U(\infty)\subset G(\infty)$. Under this correspondence, the
spherical functions of $G(\infty)$, when restricted to
$H(\infty):=\varinjlim H(n)$, coincide with the characteristic
functions of the ergodic measures.

It should be noted that all the results and constructions concerning
the ergodic measures of class ${\Cal M}$ can be translated into the
language of representation theory. In particular, weak convergence of
measures turns into uniform compact convergence of spherical
functions, and the problem that is solved in the present paper is a
particular case of the following general problem in representation
theory of ``big'' (inductive limit) groups: Given a (irreducible)
unitary representation of an inductive limit group ${\Cal
G}=\varinjlim{\Cal G}(n)$, study its approximation by (irreducible)
unitary representations of the growing subgroups ${\Cal G}(n)$ as
$n\to\infty$.

(About various asymptotic results in representation theory of ``big''
groups, see the series of papers by the authors and S. V. Kerov,
\cite{VK1--VK3, Ke, KV, O1--O7}.)

A remarkable feature of representation theory of ``big'' groups is
its connection with analytic problems of total positivity theory (see
Thoma \cite{T}, Boyer \cite{Bo}, Vershik and Kerov \cite{VK1,~VK2},
Voiculescu \cite{Vo1}; for a systematic exposition of
total positivity theory, see Karlin's fundamental monograph \cite{K}).

To describe this connection, we need the definition of totally
positive (TP) functions: these are real-valued functions
$\varphi(t)$ on ${\Bbb R}$ such that for any $n$ and any choice of
real numbers $t_1<\cdots<t_n$ and $s_1<\cdots<s_n$, the determinant
of the $n\times n$ matrix $[\varphi(t_i-s_j)]$ is nonnegative.
Similarly, one also defines two-sided and one-sided TP sequences:
these are TP functions defined on the 1-dimensional lattice ${\Bbb
Z}\subset{\Bbb R}$ or its nonnegative part ${\Bbb Z}_+$, respectively.

It turns out that the Fourier transform of TP functions $\varphi$ on
${\Bbb R}$ leads exactly to functions $F$ of the form (0.3): this
claim is Schoenberg's main theorem in \cite{S2}. Thus, there is a
correspondence $M\leftrightarrow\varphi$ between ergodic measures
$M\in{\Cal M}$ and TP functions $\varphi$ (normalized by the
condition $\int\varphi(t)\,dt=1$) that can be stated as follows:
$\varphi(t)\,dt$ coincides with the distribution of (any) diagonal
entry $B_{ii}$, where $B$ stands for a random Hermitian matrix
distributed according to $M$.

Further, there exists a similar correspondence between characters of
the group $U(\infty)$ or of the infinite symmetric group $S(\infty)$
and two-sided or one-sided TP sequences, respectively.

This correspondence can be used in two directions:

On the one hand, old theorems in total positivity theory, obtained by
analytic tools in \cite{ASW, E1, E2, S2, T, K},
can be applied to representation-theoretic
problems (classification of spherical functions or characters). Such
an approach is adopted by Thoma \cite{T}, Boyer \cite{Bo}, and
Pickrell \cite{Pi2}.

But on the other hand, if we can classify spherical functions or
characters independently, then we can prove some theorems on TP
functions or TP sequences in a new way.

It is the second approach that is adopted in Vershik and Kerov's
papers \cite{VK1, VK2}, and also in the present paper: as a
corollary of our main result, we obtain a new derivation of
Schoenberg's classification for TP functions on the real line.

Also note a recent paper by Okounkov \cite{Ok}, where a different
(direct represen\-ta\-tion-theoretic) method is used; the result of
\cite{Ok} gives yet another way to classify one-sided TP sequences.

\subhead\rom{3.}\quad The method \endsubhead
Now let us describe our techniques in more detail.

The starting point of our approach is a general approximation theorem
for ergodic measures, see Theorem 3.2 below\footnote{There is also
another version of this theorem that deals with spherical
functions, see Theorem 3.5 below.}. When applied to our concrete
situation, it implies that any ergodic measure $M\in{\Cal M}$ can be
approximated by a sequence $\{M_n\}$ of orbital measures (see Theorem
3.3 below). Then the main problem is to understand what sequences
$\{M_n\}$ are weakly convergent, as $n\to\infty$, and what are their
limits.

An attempt to do this was made in \cite{V}, but the solution proposed
there was incomplete, because of a gap in the calculations.
Nevertheless, the method itself was correct, and a refinement of the
calculations leads to the right result.

To study the weak convergence of orbital measures $M_n$, we deal with
their characteristic functions $f_n$,
$$
f_n(A)=\int_{B\in\Omega_n}e^{i\op{tr}(AB)}M_n(dB),\qquad A\in H(n),\tag0.4
$$
where $\Omega_n\subset H(n)$ is the $U(n)$-orbit carrying $M_n$. We
calculate the Taylor decomposition at the origin of a characteristic
function $f_n$ and then analyze the asymptotics of its Taylor
coefficients as $n\to\infty$. Note that these Taylor coefficients are
nothing but the moments of the measure $M_n$. (In fact, due to the
symmetry of $f_n$, it is convenient to rewrite its Taylor
decomposition as a series of Schur polynomials.) A non evident fact
is that weak convergence of orbital measures can always be controlled
by the moments. Such a phenomenon was first discovered in \cite{VK2};
there it was applied to an allied problem: classifying characters of
the group $U(\infty)$.

Note that we are able to generalize our results to the spaces of real
symmetric and quaternionic Hermitian matrices. Instead of Schur
polynomials, we must then use Jack symmetric polynomials. As a
further generalization, one could consider all matrix spaces of the
form $H=\varprojlim H(n)$, where $H(n)$ ranges over one of the 10
series of classical symmetric spaces of Euclidean type, the parameter
$n$ being the rank of the space (see \cite{O5, Pi1} and
especially \cite{Pi2} for a discussion of these spaces $H$).
Characteristic functions of orbital measures on the spaces $H(n)$ are
sometimes called generalized Bessel functions; for $3$ of the $10$ spaces
$H$, they can be expressed in terms of elementary functions, and for the
$7$ other spaces, these are certain multidimensional special functions.
The problem consists in studying their limiting behavior as
$n\to\infty$. We conjecture that our approach can be transferred to
all the spaces $H$.

\subhead\rom{4.}\quad Contents \endsubhead The present paper is organized as
follows.

In \S1, we introduce ergodic measures and explain their
relationship with spherical unitary representations.

In \S2, we formulate the so-called Multiplicativity Theorem
for characteristic functions of ergodic measures (Theorem 2.1). This
important result states that an $U(\infty)$-invariant probability
measure $M$ on $H$ is ergodic if and only if its characteristic
function is multiplicative in a certain sense; it follows that the
set of ergodic measures is closed under convolution. Then we state
the classification result (Theorem 2.9), which shows that any ergodic
measure is a convolution product of certain ``elementary'' measures.
We also discuss a number of corollaries of these theorems.

Section 3 is devoted to a proof of a general approximation theorem
that has been announced in \cite{V} (Theorem 3.2 below). This result
may be viewed as an ergodic theorem for actions of general inductive
limits of compact groups.

In \S4, we state the main result of the paper, Theorem 4.1, and
outline its proof.

Section 5 contains a preliminary result, needed for the proof of
Theorem 4.1. There we decompose the characteristic function of an
orbital measure into a series of Schur symmetric
functions.\footnote{It would be interesting to find analogues of this
result for all generalized Bessel functions.}

In \S6, we prove Theorem 4.1.

In \S7, following Pickrell's arguments in \cite{Pi2}, we
establish the equivalence of two classification problems: that of
ergodic invariant measures $M\in{\Cal M}$ and that of totally
positive functions $\varphi$ on ${\Bbb R}$. Note that our own
contribution to the results of this section is very modest. Here we
only aimed to clarify some technical details of the correspondence
$M\leftrightarrow\varphi$ and to explain how our main result implies
Schoenberg's classification.\footnote{Note, however, that we are not
completely satisfied by the way we obtain Schoenberg's theorem. We
think there should exist a more conceptual derivation of this result
by approximation methods. A similar remark can be made {\it \`a propos\/}
TP sequences as well.}

In the final \S8, we explain a connection between our main theorem
and the main result of the remarkable work \cite{CS} by Curry and
Schoenberg.

A preliminary version of the present paper appeared as preprint
\cite{OV}.

\subhead\rom{5.}\quad Acknowledgements \endsubhead
Already in the sixties, the
second named author (A.~V.) had discussions with F. A. Berezin about
the subject of the present work and its links with mathematical
physics.

The preparation of the preliminary version \cite{OV} of the paper was
completed while the first named author (G.~O.) was a guest of the
Department of Mathematical Sciences, University of Tokyo. It is a
pleasure to thank Masatoshi Noumi and the Department of Mathematical
Sciences for the kind invitation to Tokyo and hospitality.

Special thanks are due to Andre\u\i\ Okounkov
who brought to our attention the
important relation with the work \cite{CS} by Curry and Schoenberg.

The authors also acknowledge the financial support by the
International Science Foundation and by the Russian Foundation for Basic
Research during the final phase of our work.

\head\S1.\enspace Ergodic measures and spherical representations \endhead

Let $H(n)$ denote the real vector space formed by complex Hermitian
$n\times n$ 
matrices, $n=1,2,\dots$. There is a natural embedding
$$
H(n)\to H(n+1),\qquad A\mapsto\bmatrix A&0\\ 0&0\endbmatrix,
$$
and we denote by
$H(\infty)$ the corresponding inductive limit space $\varinjlim
H(n)$. Then $H(\infty)$ is identified with the space of infinite
Hermitian matrices with a finite number of nonzero entries. We
equip $H(\infty)$ with the inductive limit topology. In particular, a
function $f\:H(\infty)\to{\Bbb C}$ is continuous if its restriction to
$H(n)$ is continuous for any $n$.

Let $H$ stand for the space of all infinite Hermitian matrices. For
$A\in H$ and $n=1,2,\dots$, we denote by $\theta_n(A)\in H(n)$ the
upper left $n\times n$ corner of $A$. Using the projections
$\theta_n\:H\to H(n)$, $n=1,2,\dots$, we may identify $H$ with the
projective limit space $\varprojlim H(n)$. We equip $H$ with the
corresponding projective limit topology. In particular, $H$ is a
Borel space.

There is a natural pairing
$$
H(\infty)\times H\to{\Bbb R},\qquad(A,B)\mapsto\op{tr}(AB)\.
$$
Using it, we may regard $H$ as the algebraic dual space of $H(\infty)$.

Note that $H$ can be identified, in an obvious manner, with the space
${\Bbb R}^{\infty}={\Bbb R}\times{\Bbb R}\times\cdots$. Under this
identification, $H(\infty)\subset H$ turns into ${\Bbb
R}^{\infty}_0:=\bigcup_{n\ge1}{\Bbb R}^n$, and the pairing defined above
becomes the standard pairing between ${\Bbb R}^{\infty}_0$ and ${\Bbb
R}^{\infty}$.

Given a Borel probability measure $M$ on $H$, we define its Fourier
transform, or characteristic function, as the following function on
$H(\infty)$:
$$
f(A)=f_M(A)=\int_H e^{i\op{tr}(AB)}M(dB)\.\tag1.1
$$

We need the following statement:

\proclaim{Proposition 1.1} The Fourier transform\/ \rom{(1.1)}
establishes a bijective correspondence between Borel probability
measures $M$ on $H$ and continuous positive definite normalized
functions $f$ on $H(\infty)$.
\endproclaim

Here ``normalized'' means that $f(0)=1$.

\demo{Proof} This is an immediate corollary of Kolmogorov's
consistency theorem and Bochner's theorem. Indeed, let us identify
$(H(\infty),H)$ with $({\Bbb R}^{\infty}_0,{\Bbb R}^{\infty})$. Then
Kolmogorov's theorem implies that Borel probability measures on
${\Bbb R}^{\infty}$ are just the projective limits of probability
measures on the ${\Bbb R}^n$'s, and Bochner's theorem allows us to
restate this fact in terms of characteristic functions. \qed
\enddemo

Let $U(n)$ be the group of unitary $n\times n$ matrices,
$n=1,2,\dots$. For any $n$, we embed $U(n)$ into $U(n+1)$ using the
mapping $u\mapsto\bmatrix u&0\\ 0&1\endbmatrix$. Let
$U(\infty)=\varinjlim U(n)$ denote the corresponding inductive limit
group. We regard $U(\infty)$ as the group of infinite unitary
matrices $u=[u_{ij}]_{i,j=1}^{\infty}$ with a finite number of
entries $u_{ij}\ne\delta_{ij}$. The group $U(\infty)$ acts by
conjugations both on $H(\infty)$ and $H$, and the pairing between
these two spaces is clearly $U(\infty)$-invariant.

Let us recall some basic facts concerning ergodic measures.

\definition{Definition 1.2} Let ${\Cal X}$ be a Borel space, let ${\Cal
G}$ be a group of Borel transformations of ${\Cal X}$, and let $M$ be
a ${\Cal G}$-invariant probability Borel measure on ${\Cal X}$.
\roster
\item"(i)" A Borel subset ${\Cal Y}\subset{\Cal X}$ is said to be {\it
${\Cal G}$-invariant\/} mod $0$ if $M((g{\Cal Y}\)\triangle{\Cal Y})=0$
for any $g\in{\Cal G}$.
\item"(ii)" The measure $M$ is said to be {\it ergodic\/} with respect to
${\Cal G}$ if the $M$-volume of any ${\Cal G}$-invariant mod $0$ Borel
subset ${\Cal Y}\subset{\Cal X}$ is equal either to $0$ or $1$.
\endroster
\enddefinition

\proclaim{Proposition 1.3} The following conditions on a Borel
probability measure $M$ on ${\Cal X}$ are equivalent\/\rom:
\roster
\item"(i)" $M$ is ergodic\/\rom;
\item"(ii)" $M$ is an extreme point of the convex set formed by all
Borel probability measures on ${\Cal X}$\rom;
\item"(iii)" the subspace of ${\Cal G}$-invariant vectors in the
Hilbert space $L^2({\Cal X},M)$ is exhausted by constant functions.
\endroster
\endproclaim

\demo{Proof} See, e.g., Phelps \cite{Ph, Proposition 10.4}.\qed
\enddemo

For compact group actions on locally compact spaces, ergodic measures
coincide with invariant measures supported by the orbits, so that for
such actions, classifying orbits and classifying ergodic measures are
equivalent problems. However, for ``big'' groups like $U(\infty)$,
these two problems are quite different: the first one seems to be out
of reach, while the second one has, in certain cases, a nice
solution. According to a general principle of ergodic theory, ergodic
measures may be viewed as true substitutes of orbits.

Let ${\Cal M}$ denote the set of all ergodic $U(\infty)$-invariant
Borel probability measures $M$ on $H$, and let ${\Cal F}$ denote the
set of the corresponding characteristic functions $f=f_M$.

The following claim, which follows from Propositions 1.1 and 1.3,
gives a useful characterization of the class ${\Cal F}$:

\proclaim{Proposition 1.4} The functions $f\in{\Cal F}$ are exactly
the extreme points of the convex set formed by all continuous
$U(\infty)$-invariant positive definite normalized functions on
$H(\infty)$. \qed
\endproclaim

Let $G(n)=U(n)\ltimes H(n)$ be the semidirect product of $U(n)$ and
the additive group of the vector space $H(n)$. The elements $g\in
G(n)$ are the pairs $(u,A)\in U(n)\times H(n)$ with the
multiplication rule
$$
(u,A)\cdot(v,B)=(uv,v^{-1}Av+B)\.
$$
In the same way, we define the group $G(\infty)=U(\infty)\ltimes
H(\infty)$. The group $G(\infty)$ can be also viewed as the inductive
limit group $\varinjlim G(n)$, and we equip it with the inductive
limit topology. Using the embeddings
$$
u\mapsto(u,0)\quad\text{and}\quad A\mapsto(1,A),\qquad u\in U(\infty),
\quad A\in H(\infty),
$$
we may identify $U(\infty)$ and $H(\infty)$ with the corresponding
subgroups in $G(\infty)$.

Let $T$ be a unitary representation of $G(\infty)$ in a Hilbert space
${\Cal H}(T)$ (we tacitly assume $T$ is continuous with respect to
the inductive limit topology of the group $G(\infty)$). Then $T$ is
said to be {\it spherical\/} if it is irreducible and the subspace
${\Cal H}(T)^{U(\infty)}$ of $U(\infty)$-invariants in ${\Cal H}(T)$
is nonzero.

Suppose $T$ is spherical. It is then known that the space ${\Cal
H}(T)^{U(\infty)}$ is one-dimensional (see Olshanski \cite{O5,
Theorem 23.6}). A vector $h\in{\Cal H}(T)^{U(\infty)}$ of norm $1$ is
called a {\it spherical vector\/} of $T$, and the corresponding
matrix element
$$
\varphi_T(g)=(T(g\)h,h),\qquad g\in G(\infty),
$$
is called the {\it spherical function\/} of $T$. Since $\varphi_T$
does not change when $h$ is multiplied by a complex number of
absolute value $1$, $\varphi_T$ is an invariant of $T$; moreover, it
uniquely determines $T$. Further, as $\varphi_T$ is bi-invariant
with respect to the subgroup $U(\infty)\subset G(\infty)$, it is
uniquely determined by its restriction $\varphi_T|H(\infty)$ to
$H(\infty)$, which is a continuous $U(\infty)$-invariant function.

\proclaim{Proposition 1.5} There is a natural bijective
correspondence $T\leftrightarrow M$ between \rom(equivalence classes
of\rom) spherical representations $T$ of the group $G(\infty)$ and
ergodic measures $M\in{\Cal M}$. Under that correspondence, the
functions $\varphi_T|H(\infty)$ coincide with the characteristic
functions $f_M$. Given $M\in{\Cal M}$, the corresponding
representation $T$ can be realized in the Hilbert space $L^2(H,M)$ as
follows\/\rom:
$$
\alignat2
(T(u\)\Psi)(B)&=\Psi(u^{-1}Bu),&\qquad u&\in U(\infty)\subset G(\infty),\\
(T(A\)\Psi)(B)&=e^{i\op{tr}(AB)}\Psi(B),&\qquad A&\in
H(\infty)\subset G(\infty),
\endalignat
$$
where $\Psi\in L^2(H,M)$ and $B\in H$. Finally, in this realization,
the spherical vector is the constant function $\Psi_0(B)\equiv1$.
\endproclaim

\demo{Proof} Indeed, it is enough to remark that the
spherical functions, when restricted to $H(\infty)$, are just the
extreme $U(\infty)$-invariant continuous positive definite functions,
normalized at $0\in H(\infty)$, and then apply Proposition 1.4.\qed
\enddemo

\head\S2.\enspace The Classification Theorem \endhead

Denote by $D(\infty)$ the subspace of diagonal matrices in
$H(\infty)$. An element of $D(\infty)$ will be written as
$\op{diag}(a_1,a_2,\dots)$, where $a_1,a_2,\dots\in{\Bbb R}$ and
$a_k=0$ for $k$ large enough. Since any matrix in $H(\infty)$ can be
diagonalized under the action of the group $U(\infty)$, a
$U(\infty)$-invariant function $f$ on $H(\infty)$ is uniquely
determined by $f\mi D(\infty)$, the restriction of $f$ to $D(\infty)$.

Assume $f$ is an $U(\infty)$-invariant function on $H(\infty)$,
$f(0)=1$. Let us say that $f$ is {\it multiplicative\/} if
$$
f(\op{diag}\(a_1,a_2,\dots))=F(a_1)F(a_2)\cdots,\tag2.1
$$
where $F$ is a function on ${\Bbb R}$ such that $F(0)=1$. In other
words,
$$
f(A)=F(a_1)F(a_2)\cdots,\qquad A\in H(\infty),\tag2.2
$$
where $a_1,a_2,\dots$ are the eigenvalues of $A$. Note that
$F(a)=f(\op{diag}\(a,0,0,\dots))$, i.e., $F=f\mi H(1)$.

\proclaim{Theorem 2.1 \rm(Multiplicativity Theorem)} Let $f$ be a
continuous $U(\infty)$-invariant positive definite function on
$H(\infty)$, $f(0)=1$. Then $f$ is extreme\/ \rom(i.e., $f\in{\Cal
F}$\rom) if and only if $f$ is multiplicative.
\endproclaim

\demo{Proof} See, e.g., Olshanski \cite{O5, Theorem 23.8}. There the
Multiplicativity Theorem is stated for allied groups (like
$GL(n,{\Bbb C})$) but the proof is exactly the same. \qed
\enddemo

Not that there are many theorems of this type and a lot of
different methods to prove them: see Thoma \cite{T}, Ismagilov
\cite{I1} and \cite{I2} (Ismagilov's method is also explained in
Olshanski \cite{O4, section 2.5}), Nessonov \cite{N}, Olshanski
\cite{O1}, Pickrell \cite{Pi1}, Vershik--Kerov \cite{VK4}, Voiculescu
\cite{Vo1} and \cite{Vo2}, Stratila--Voiculescu \cite{SV}.

Theorem 2.1 implies that a function $f\in{\Cal F}$ is uniquely
determined by the function $F:=f\mi H(1)$, which is a continuous
function in one real variable. Note that
$$
F(a)=F_M(a)=\int_H e^{iaB_{11}}M(dB),\qquad a\in{\Bbb R},\;B\in H,\tag2.3
$$
where $M\in{\Cal M}$ corresponds to $f$.

\definition{Definition 2.2} Let ${\Cal F}_1=\{F\}$ denote the class of
all functions on ${\Bbb R}$ of the form $F=f\mi H(1)$ where $f$ ranges
over ${\Cal F}$. In other words,
a function $F$ on ${\Bbb R}$ belongs to ${\Cal
F}_1$ if and only if $F$ is continuous, $F(0)=1$, and the
corresponding function (2.2) on $H(\infty)$ is positive definite.
\enddefinition

Clearly, the classification of the measures $M\in{\Cal M}$ is reduced
to that of the functions $F\in{\Cal F}_1$.

Theorem 2.1 has a number of corollaries.

\proclaim{Corollary 2.3} The class ${\Cal F}_1$ is stable under
pointwise multiplication.\qed
\endproclaim

This implies that the class ${\Cal F}$ is also stable under
multiplication, and the class ${\Cal M}$ is stable under convolution.

\proclaim{Corollary 2.4} If a sequence $F_1,F_2,\dots\in{\Cal F}_1$
pointwise converges to a continuous function $F$ on ${\Bbb R}$, then
$F\in{\Cal F}_1$.\qed
\endproclaim

\proclaim{Corollary 2.5} For a real $\gamma$, the Dirac measure
concentrated at $\gamma\cdot1\in H$ belongs to ${\Cal M}$. The
corresponding function $F\in{\Cal F}_1$ is
$$
F(a)=e^{i\gamma a},\qquad a\in{\Bbb R}\.\tag2.4
$$
\endproclaim

\demo{Proof} Since scalar matrices are $U(\infty)$-invariant, our
measure is invariant. It is clearly ergodic. It is evident that the
corresponding function from ${\Cal F}_1$ is given by (2.4). \qed
\enddemo

\proclaim{Corollary 2.6} Given $\gamma\ge0$, let $M$ be the Gaussian
distribution on $H$ such that, for a matrix $B\in H$, the diagonal
entries $B_{ii}$ and the off-diagonal entries
$\op{Re}B_{ij}$, $\op{Im}B_{ij}$, $i<j$, are
independent Gaussian variables with mean $0$ and variance $\gamma$.
Then $M\in{\Cal M}$, and the corresponding function $F\in{\Cal F}_1$ is
$$
F(a)=e^{-\gamma a^2\!/2},\qquad a\in{\Bbb R}\.\tag2.5
$$
\endproclaim

\demo{Proof} It is easily verified that $M$ is invariant and its
characteristic function satisfies (2.1), where $F$ is given by
(2.5).\qed
\enddemo

\proclaim{Corollary 2.7} Let $\omega$ denote the Gaussian measure on
${\Bbb C}$ with density given by $\pi^{-1}\exp\(-|z|^2)$, $z\in{\Bbb
C}$. Given $y\in{\Bbb R}$, let $M$ be the image of the Gaussian
product measure $\omega^{\otimes\infty}$ on ${\Bbb C}^{\infty}$ under
the following Borel mapping ${\Bbb C}^{\infty}\to H$\rom:
$$
{\Bbb C}^{\infty}\ni\xi\mapsto y(-1+\xi^*\xi)=B\in H\.\tag2.6
$$
\rom(Here $\xi=(\xi_1,\xi_2,\dots)\in{\Bbb C}^{\infty}$ is
regarded as a row vector, so that the $(i,j)$-entry of the matrix
$B$ is equal to $y(-1+{\ov\xi_i}\xi_j)$, $i,j=1,2,\dots$.\rom)
Then $M\in{\Cal M}$, and the corresponding function $F\in{\Cal F}_1$
is given by
$$
F(a)=\frac{e^{-iya}}{1-iya},\qquad a\in{\Bbb R}\.\tag2.7
$$
\endproclaim

\demo{Proof} The invariance of $M$ is obvious. A direct calculation
shows that the characteristic function of $M$ satisfies (2.1) with
$F$ given by (2.7). \qed
\enddemo

We shall call the measures defined in Corollaries 2.5--2.7 the {\it
elementary ergodic measures.}

Note that we could omit the term $-1$ in (2.6) and then the numerator
in (2.7) would disappear. However, due to this term, the mean of $M$
is equal to zero and the function (2.7) has the property
$$
F(a)=\frac{e^{-iya}}{1-iya}=1-\frac32\,a^2y^2+O(y^3)\quad
\text{as }\,y\to0,\tag2.8
$$
i.e., the term of degree $1$ in $y$ vanishes. This is important for the
following construction.

\proclaim{Proposition 2.8} The class ${\Cal F}_1$ contains all
functions of the form
$$
F_{\gamma_1,\gamma_2,x}(a)=e^{i\gamma_1a-\gamma_2a^2\!/2}
\prod_k\frac{e^{-ix_ka}}{1-ix_ka},\tag2.9
$$
where $\gamma_1\in{\Bbb R}$ and $\gamma_2\ge0$ are arbitrary
constants, and $x=(x_1,x_2,\dots)$ is a sequence of real numbers
such that $\sum x_k^2<\infty$.
\endproclaim

\demo{Proof} Suppose first the sequence $x$ is finite. Then (2.9) is
a finite product of functions that belong to the class ${\Cal F}_1$
due to Corollaries 2.5, 2.6, and 2.7. By Corollary 2.2, their product
belongs to ${\Cal F}_1$ also. Finally, if the sequence $x$ is
infinite, then the assumption $\sum x_k^2<\infty$ and the estimate
(2.8) imply that the product in (2.9) is convergent for any
$a\in{\Bbb R}$; moreover, the result is a continuous function. Then
Corollary 2.4 implies that this function belongs to ${\Cal F}_1$. \qed
\enddemo

\example{Comments} 1) The order of the $x_k$'s is unessential, so
that $x=(x_k)$ is rather a point configuration (or multiset) than a
sequence.

2) If $\sum|x_k|<\infty$ then (2.9) may be rewritten as
$$
F_{\gamma_1,\gamma_2,x}(a)=e^{i{\ov\gamma_1}a-\gamma_2a^2\!/2}
\prod_k\frac1{1-ix_ka},\qquad \ov\gamma_1:=\gamma_1-\sum_k x_k\.\tag2.10
$$

3) The function (2.9) admits a holomorphic continuation to the
horizontal strip $\{z\in{\Bbb
C}\mid|\op{Im}z|<\varepsilon\}$, where $\varepsilon^{-1}=\sup|x_k|$.

4) The parameters $\gamma_1$, $\gamma_2$, $x$ are uniquely determined
by the function $F_{\gamma_1,\gamma_2,x}$.

5) The characteristic function (2.2) corresponding to the function
(2.9) can be written as
$$
f(A)=e^{i\gamma_1\op{tr}A-\gamma_2\op{tr}(A^2)/2}
\det\bigg(\prod_k\frac{e^{-ix_kA}}{1-ix_kA}\bigg),
\qquad A\in H(\infty)\.\tag2.11
$$
\endexample

\proclaim{Theorem 2.9 \rm(Classification Theorem)} The class
${\Cal F}_1$ is exhausted by the functions of the form\/ \rom{(2.9)}.
Thus, the characteristic functions of the ergodic
$U(\infty)$-invariant Borel probability measures on $H$ are just the
functions of the form\/ \rom{(2.11)}, where $\gamma_1,\gamma_2,
x_1,x_2,\dots$ are real parameters such that $\gamma_2\ge0$ and
$\sum x_k^2<\infty$.
\endproclaim

This result gives a description of ergodic measures $M\in{\Cal M}$:
any such $M$ is a convolution of the elementary ergodic measures
constructed in Corollaries 2.5--2.7. The proof will be given in \S4:
we shall derive Theorem 2.9 from a more general result, Theorem 4.1.

Note that the elementary ergodic measures of Corollaries 2.5 and 2.6
are infinitely divisible with respect to convolution, whereas those of
Corollary 2.7 are not.

\remark{Remark\/ \rm2.10} The construction of Corollary 2.7 can be
generalized as follows. Fix $k=1,2,\dots$, consider the space ${\Bbb
C}^{k\times\infty}$ of all complex matrices $\Xi$ with $k$ rows and
infinitely many columns, and equip it with the Gaussian product
measure
$$
\omega^{k\times\infty}:=\omega^{\otimes\infty}\otimes\cdots
\otimes\omega^{\otimes\infty}.\tag2.12
$$
Let $z,x_1,\dots,x_k$ be any real numbers, and let
$X=\op{diag}\(x_1,\dots,x_k)$ denote the diagonal matrix of order $k$
with diagonal entries $x_1,\dots,x_k$. Consider the mapping
$$
{\Bbb C}^{k\times\infty}\ni\Xi\mapsto z\cdot1+\Xi^*X\Xi=B\in H\.\tag2.13
$$
Then the image of the Gaussian measure (2.12) under the mapping
(2.13) is an ergodic measure $M_{z;x_1,\dots,x_k}$. Its parameters
are
$$
\gamma_1=z+x_1+\cdots+x_k,\quad\gamma_2=0,\quad x=(x_1,\dots,x_k,
0,0,\dots)\.\tag2.14
$$
$M_{z;x_1,\dots,x_k}$, $k=1,2,\dots$, form a weakly dense subset of
${\Cal M}$. Further, let $H_{\le k}$ denote the closed subspace of
$H$ formed by matrices of rank $\le k$. One can prove that the
measures $M_{0;x_1,\dots,x_k}$ are just those measures $M\in{\Cal
M}$ that are supported by $H_{\le k}$. Finally, note that the
measures $M_{0;x_1,\dots,x_k}$ with $x_1=\dots=x_k$ are
infinite-dimensional analogs of the well-known Wishart distributions
(see e.g., Muirhead \cite{Mu}).
\endremark

\remark{Remark\/ \rm2.11} Let $H_+(n)\subset H(n)$ and $H_+\subset H$
denote the subsets of nonnegative definite matrices. Then $H_+$
coincides with $\varprojlim H_+(n)$ and is a closed cone in $H$. One
can prove that a measure $M\in{\Cal M}$ is supported by $H_+$ if and
only if its parameters satisfy the conditions
$$
\gamma_2=0,\qquad x_1\ge0,\quad x_2\ge0,\quad\dots,\quad
\sum x_k\le\gamma_1<\infty\.\tag2.17
$$
\endremark

\remark{Remark\/ \rm2.12} Let $M\in{\Cal M}$ and let $F\in{\Cal F}_1$ be
the corresponding function (2.9). We may regard $(H,M)$ as a
probability space and the matrix elements $B_{ij}$ of a matrix $B\in
H$ as random variables. Let $\mu$ denote the distribution of the real
random variable $B_{11}$; then $F$ is the characteristic function of
$\mu$. We know that $M$ is completely determined by $\mu$. Let us
describe the structure of $\mu$. It follows from (2.9) that $\mu$ is
the convolution of a (not necessarily centered) normal distribution
with a family of distributions possessing characteristic functions of the
form (2.7). It is easily verified that if $y>0$, then the distribution
with characteristic function (2.7) is supported by the half-line
$t\ge-y$ and has density
$$
t\mapsto y^{-1}e^{-y^{-1}(t+y)},\qquad t\ge-y\.\tag2.18
$$
This is the shifted exponential distribution with variance $y^2$ and
mean 0. If the parameter $y$ is negative, it suffices to replace $t$
by $-t$. Thus, $\mu$ is the convolution of a normal distribution with
a family of modified exponential distributions.
\endremark

\remark{Remark\/ \rm2.13} Let $M\in{\Cal M}$ be an ergodic measure for which
not all the parameters $x_k$ vanish. Then the characteristic
function $f(A)$, see (2.11), cannot be factorized into a product of
factors each of which depends on a single matrix element $A_{ij}$
only. This means that $M$ is not a product measure, so that the
matrix elements $B_{pq}$, as random variables defined on $(H,M)$, are
not independent on the whole (in the particular case of the measures
$M_{z;x_1,\dots,x_k}$, this is seen from their construction via the
mapping (2.13)). However, certain matrix elements are independent.
For instance, it is evident that the diagonal elements are
independent. More generally, for any $(p_1,q_1),\dots,(p_m,q_m)$
such that $i\ne j$ implies $p_i\ne p_j$, $q_i\ne q_j$, $p_i\ne q_j$,
the matrix elements $B_{p_1q_1},\dots,B_{p_mq_m}$ are independent.
This claim can be deduced from (2.2).
\endremark
\medskip

According to Proposition 1.5, Theorem 2.9 also gives a complete
description of spherical representations of the group $G(\infty)$. We
shall now discuss the possibility of extending spherical
representations to some topological completions of the inductive
group $G(\infty)$.

Let us regard $U(\infty)$ as a group of unitary operators in the
complex coordinate Hilbert space $\ell_2$, and let
$\ov{U}(\infty)\supset U(\infty)$ stand for the group of all
unitary operators in $\ell_2$. Note that $\ov{U}(\infty)$ is a
topological group with respect to the weak operator topology (which
coincides on unitary operators with the strong operator topology).
Further, let $H(\infty)_1\supset H(\infty)$ (respectively,
$H(\infty)_2\supset H(\infty)$) denote the space of the trace class
(respectively, Hilbert--Schmidt) Hermitian operators in $\ell_2$,
equipped with the topology defined by the trace norm
$\|\bdot\|_1$ (respectively, by the Hilbert--Schmidt norm
$\|\bdot\|_2$).

One can check that the actions
$$
\ov{U}(\infty)\times H(\infty)_1\to H(\infty)_1,\qquad
\ov{U}(\infty)\times H(\infty)_2\to H(\infty)_2,
$$
where $(u,A)\mapsto uAu^{-1}$, are continuous (cf.~Shale \cite{Sha}).
Thus we may form the semidirect products
$$
G(\infty)_1=\ov{U}(\infty)\ltimes H(\infty)_1,
\qquad G(\infty)_2=\ov{U}(\infty)\ltimes H(\infty)_2,
$$
which are topological groups with respect to the corresponding
product topologies. Note that the inductive limit group
$G(\infty)=U(\infty)\ltimes H(\infty)$ is contained as a dense
subgroup both in $G(\infty)_1$ and $G(\infty)_2$.

\proclaim{Corollary 2.14} {\rm (i)} Any spherical representation of
the group $G(\infty)$ admits an extension to a continuous
representation of the topological group $G(\infty)_1$.

{\rm (ii)} It can be continued to $G(\infty)_2$ if and only if the
parameter $\gamma_1$ in\/ \rom{(2.11)} vanishes.
\endproclaim

\demo{Proof} (i) Indeed, for any function $F$ of the form (2.9), the
corresponding function $f$ on $H(\infty)$, which is defined by (2.11),
can be extended to a continuous function on $H(\infty)_1$. Then the
latter function can be extended to a continuous
$\ov{U}(\infty)$-bi-invariant function on the group
$G(\infty)_1$, and our claim follows.

(ii) We argue as in (i), using the fact that $f$ is continuous with
respect to the Hilbert--Schmidt norm if and only if $\gamma_1=0$. \qed
\enddemo
\head\S3.\enspace Approximation of ergodic measures by orbital measures
\endhead

Let ${\Cal X}$ be a separable metric space and $C({\Cal X})$ be the
Banach space of bounded continuous functions on ${\Cal X}$. By a {\it
measure\/} on ${\Cal X}$ we shall always mean a Borel probability
measure. Recall that a sequence $\{\nu_n\}$ of measures on ${\Cal X}$
is said to be {\it weakly convergent to a measure\/} $\nu$ (notation:
$\nu_n\Rightarrow\nu$) if
$\langle\psi,\nu_n\rangle\to\langle\psi,\nu\rangle$ as $n\to\infty$
for any $\psi\in C({\Cal X})$, see, e.g., \cite{Bi, Pa1}.

\proclaim{Proposition 3.1} There exists a countable set $\Psi\subset
C({\Cal X})$ of functions with the following property: a sequence of
measures $\nu_1,\nu_2,\dots$ weakly converges to a measure $\nu$
provided convergence occurs for any function $\psi\in\Psi$.
\endproclaim

\demo{Proof} See, e.g., \cite{Pa1, Chapter II, Theorem 6.6}. \qed
\enddemo

Let ${\Cal K}(1)\subset{\Cal K}(2)\subset\dots$ be an ascending
chain of compact groups, and let ${\Cal K}=\varinjlim{\Cal K}(n)$ be
the corresponding inductive limit topological group. We assume there
is a jointly continuous action $(u,x)\mapsto u\cdot x$ of the group
${\Cal K}$ on the space ${\Cal X}$. Let $m_n$ denote the normalized
Haar measure on ${\Cal K}(n)$. Given a point $x\in{\Cal X}$, let
$m_n(x)$ denote the image of $m_n$ under the mapping $u\mapsto u\cdot
x$, where $u$ ranges over ${\Cal K}(n)$, i.e., $m_n(x)$ is the unique
${\Cal K}(n)$-invariant probability measure supported by the orbit
${\Cal K}(n)\cdot x$.

\proclaim{Theorem 3.2 \rm(Vershik \cite{V, Theorem 1})} Let $\nu$
be an ergodic ${\Cal K}$-invariant Borel probability measure on
${\Cal X}$. Then there exists a point $x\in{\Cal X}$ such that
$m_n(x)\Rightarrow\nu$ as $n\to\infty$. Moreover, the set of all
points $x\in{\Cal X}$ with this property is of full measure with
respect to $\nu$.
\endproclaim

\demo{Proof} Given $\psi\in C({\Cal X})$, put
$$
\psi_n(x)=\langle\psi,m_n(x)\rangle=\int_{{\Cal K}(n)}\psi(u\cdot x\)
m_n(du),\qquad x\in{\Cal X},\tag3.1
$$
and let $\ov{\psi}$ denote the constant function
$$
\ov{\psi}(x)\equiv\langle\psi,\nu\rangle\cdot1=
\bigg(\int_{{\Cal X}}\psi(x\)\nu(dx)\bigg)\cdot1\.\tag3.2
$$

To prove the theorem, it is enough to check that $\psi_n\to\ov\psi$
almost everywhere for any $\psi\in C({\Cal X})$. 
Indeed, it will follow that, given an arbitrary
countable family of functions $\Psi\subset C({\Cal X})$, there exists
a subset ${\Cal X}'$ of full measure such that
$$
\langle\psi,m_n(x)\rangle\to\langle\psi,\nu\rangle\qquad\text{for any
$\psi\in\Psi$ and any $x\in{\Cal X}'$}.\tag3.3
$$
Taking as $\Psi$ a subset from Proposition 3.1, we shall obtain
$m_n(x)\Rightarrow\nu$ for all $x\in{\Cal X}'$.

Now we proceed to check that, given $\psi\in C({\Cal X})$,
we have $\psi_n\to\ov\psi$ almost everywhere.

Remark that the $\psi_n$ are bounded Borel functions. Indeed, given
$x\in{\Cal X}$ and $n$, the function $u\mapsto\psi(u\cdot x)$ is
continuous on ${\Cal K}(n)$, so that the integral in (3.1) can be
approximated by its Riemannian integral sums. It follows that
$\psi_n$ is a pointwise limit of continuous functions, so it is a
Borel function. Its boundedness is immediate.

Now it is enough to verify the following two claims:

{\it Claim\/} 1. $\psi_n\to\ov{\psi}$ in the metric of
$L^2({\Cal X},\nu)$.

{\it Claim\/} 2. For almost all points $x\in{\Cal X}$, the
following limit exists
$$
\psi_{\infty}(x):=\lim_{n\to\infty}\psi_n(x)\.
$$

Indeed, since the $\psi_n$'s are uniformly bounded, Claim 2 will
imply that $\psi_n\to\psi_{\infty}$ in $L^2({\Cal X},\nu)$, whence
$\psi_{\infty}=\ov\psi$ almost everywhere, so that
$\psi_n\to\ov\psi$ almost everywhere.

To verify Claim 1, consider the natural unitary representation of the
group ${\Cal K}$ in the Hilbert space $L^2({\Cal X},\nu)$. Let $P_n$
denote the orthoprojection in this space onto the subspace of ${\Cal
K}(n)$-invariant vectors, $n=1,2,\dots$. Then $P_n$ strongly
converges, as $n\to\infty$, to $P$, the orthoprojection onto the
subspace of ${\Cal K}$-invariant vectors. In particular, $P_n\psi\to
P\psi$ in the metric of $L^2({\Cal X},\nu)$. Since the measure $\nu$
is ergodic, the only ${\Cal K}$-invariant vectors are the constants
(this is the only place when we use the assumption that $\nu$ is
ergodic). Therefore, $P\psi$ is a constant function, which clearly
equals $\ov{\psi}$. Finally, $P_n\psi=\psi_n$, so that
$\psi_n\to\ov{\psi}$ in $L^2({\Cal X},\nu)$.

Now we shall prove Claim 2 by a method similar to the one used in the
proof of Birkhoff--Khinchine's individual ergodic theorem (see, e.g.,
Parthasarathy \cite{Pa2, \S49}).

Without loss of generality one may assume that $\psi$ is
real-valued. For $N=1,2,\dots$, put
$$
\align
E_N=E_N(\psi)&=\Big\{x\in{\Cal X}\mi
\sup_{1\le n\le N}\psi_n(x)>0\Big\},\tag3.4\\
E_\infty=E_\infty(\psi)&=\Big\{x\in{\Cal X}\mi
\sup_{1\le n\le\infty}\psi_n(x)>0\Big\}
=\bigcup_{N=1}^{\infty}E_N(\psi)\.\tag3.5
\endalign
$$
Note that each $E_N$ is a Borel subset, whence $E_\infty$ also is a
Borel subset.

The following claim is an analog of the maximal ergodic theorem.

{\it Claim\/} 3. We have
$$
\int_{E_\infty}\psi(x\)\nu(dx)\ge0\.\tag3.6
$$

Indeed, since $(E_N)$ is a monotone family of sets, it is enough to
check that
$$
\int_{E_N}\psi(x\)\nu(dx)\ge0,\qquad N=1,2,\dots.\tag3.7
$$
We can write each $E_N$ as a disjoint union of Borel subsets,
$$
E_N=E_{1N}\cup\cdots\cup E_{NN},\tag3.8
$$
where
$$
E_{mN}=\{x\in{\Cal X}\mid\psi_m(x)>0,\;\psi_i(x)\le0\,\text{ for }\,
m+1\le i\le N\}\.\tag3.9
$$
Then it is enough to show that
$$
\int_{E_{mN}}\psi(x)\nu(dx)\ge0,\quad1\le m\le N.\tag3.10
$$
To do this, we remark that $E_{mN}$ is invariant relative to the
action of ${\Cal K}(m)$. Since $\nu$ is an invariant measure, it
follows that
$$
\int_{E_{mN}}\psi(x\)\nu(dx)=\int_{E_{mN}}\psi(u\cdot x\)\nu(dx),
\qquad u\in{\Cal K}(m)\.\tag3.11
$$
Since the function $(u,x)\to\psi(u\cdot x)$ is continuous on ${\Cal
K}(m)\times E_{mN}$, we may integrate the right-hand side of (3.11)
over ${\Cal K}(m)$ (with respect to the Haar measure) and then
interchange the integrals over ${\Cal K}(m)$ and over $E_{mN}$. This
yields
$$
\int_{E_{mN}}\psi(x\)\nu(dx)=\int_{E_{mN}}\psi_m(x\)\nu(dx)\.\tag3.12
$$
By definition (3.9), $\psi_m$ is positive on $E_{mN}$, so that (3.12)
is nonnegative.

Thus, we have checked Claim 3.

Further, for arbitrary real $a<b$, put
$$
{\Cal X}_{ab}=\{x\in{\Cal X}\mid\varliminf\psi_n(x)<a<b<\varlimsup
\psi_n(x)\}\.\tag3.13
$$
This is a Borel subset. Let us establish the double inequality
$$
a\,\nu({\Cal X}_{ab})\ge\int_{{\Cal X}_{ab}}\psi(x)\,\nu(dx)\ge b\,
\nu({\Cal X}_{ab}),\tag3.14
$$
which will imply $\nu({\Cal X}_{ab})=0$.

Indeed, ${\Cal X}_{ab}$ is a ${\Cal K}$-invariant Borel subset of
${\Cal X}$. Let us replace ${\Cal X}$ by ${\Cal X}_{ab}$ and apply
Claim 3 to the functions
$$
{\psi}':=(\psi-b)|_{{\Cal X}_{ab}}\quad\text{and}\quad
{\psi}'':=(a-\psi)|_{{\Cal X}_{ab}}\.\tag3.15
$$
By definition (3.13) of ${\Cal X}_{ab}$,
$$
E_\infty({\psi}')=E_\infty({\psi}'')={\Cal X}_{ab},\tag3.16
$$
whence
$$
\int_{{\Cal X}_{ab}}{\psi}'(x)\nu(dx)\ge0,\qquad
\int_{{\Cal X}_{ab}}{\psi}''(x)\nu(dx)\ge0,\tag3.17
$$
which is equivalent to (3.14).

Finally, applying (3.14) to various couples of rational numbers
$a<b$, we see that
$$
\varliminf\psi_n(x)=\varlimsup\psi_n(x)\tag3.18
$$
almost everywhere. This completes the proof of Claim 2 and of the
theorem. \qed
\enddemo

We shall apply Theorem 3.2 to ${\Cal X}=H$ and ${\Cal K}=U(\infty)$.
Note that the assumptions of Theorem 3.2 are satisfied. Indeed, $H$
is homeomorphic to a separable metric space, because it is
essentially a copy of ${\Bbb R}^{\infty}$, and, further, the action
$U(\infty)\times H\to H$ is jointly continuous.

For $n=1,2,\dots$, denote by ${\Cal M}_n$ the set of
$U(n)$-invariant probability measures that are supported by the
$U(n)$-orbits in the space $H(n)$. These measures will be called {\it
orbital measures\/}. Since $H(n)$ is contained in $H$, we may view
orbital measures as measures on the space $H$.

\proclaim{Theorem 3.3} For any ergodic measure $M\in{\Cal M}$, there
exists a sequence $\{M_n\in{\Cal M}_n\}$ of orbital measures such
that $M_n\Rightarrow M$ as $n\to\infty$.
\endproclaim

\demo{Proof} We shall write $m_n(B)$ instead of $m_n(x)$; here $B$ is
a matrix from $H$. By the first claim of Theorem 3.2, there exists
$B\in H$ such that $m_n(B)\Rightarrow M$ as $n\to\infty$. Consider
the projections $\theta_n\:H\to H(n)$ defined at the beginning of
\S1 and remark that $\theta_n(m_n(B))$ coincides with the
orbital measure in $H(n)$ corresponding to the matrix $\theta_n(B)$.
Let us take this orbital measure as $M_n$. If $k\le n$, then $m_n(B)$
and $M_n$ have the same image under $\theta_k$; fixing $k$ and
letting $n\to\infty$, we see that the measures $\theta_k(M_n)$ on
$H(k)$ weakly converge to the measure $\theta_k(M)$. Now identify $H$
with ${\Bbb R}^{\infty}$ and recall that on ${\Bbb R}^{\infty}$ weak
convergence of probability measures is equivalent to weak convergence
of their finite-dimensional projections (see, e.g., \cite{Bi,
chapter 1, \S 3}). Applying this to our sequence $\{M_n\}$, we see
that $M_n\Rightarrow M$.\qed
\enddemo

It is convenient to analyze weak convergence of measures in terms of
characteristic functions. Assume $\nu,\nu_1,\nu_2,\dots$ are Borel
probability measures on $H$ and $f,f_1,f_2,\dots$ denote their
characteristic functions. We shall need the following simple claim,
which again is essentially a well-known fact about the space ${\Bbb
R}^{\infty}$.

\proclaim{Proposition 3.4} Weak convergence $\nu_n\Rightarrow\nu$ on
$H$ is equivalent to uniform convergence $f_n\to f$ on compact
subsets in $H(\infty)$.\qed
\endproclaim

(Note that any compact subset in $H(\infty)$ is always contained in
$H(n)$ for sufficiently large $n$.)

\demo{Proof} We again pass to finite-dimensional projections and
then use the fact that weak convergence of probability measures on
${\Bbb R}^n$ is equivalent to uniform convergence, on compact sets,
of their characteristic functions. \qed
\enddemo

To obtain Theorem 3.3, we could use, instead of Theorem 3.2, another
general result, where we have to specialize ${\Cal G}=G(\infty)$,
${\Cal K}=U(\infty)$:

\proclaim{Theorem 3.5 \rm(Olshanski \cite{O3, Theorem 2.5} and
\cite{O5, Theorem 22.10})} Let ${\Cal G}=\varinjlim{\Cal G}(n)$ be
an inductive limit of separable locally compact groups and ${\Cal K}$
be a subgroup of ${\Cal G}$. Let ${\Cal K}(n)={\Cal K}\cap{\Cal
G}(n)$, so that ${\Cal K}=\varinjlim{\Cal K}(n)$.

Then any extreme ${\Cal K}$-bi-invariant continuous positive definite
function $f$ on ${\Cal G}$, normalized at unity, can be
approached, uniformly on compact sets of the group ${\Cal G}$, by a
sequence $\{f_n\}$ of extreme ${\Cal K}(n)$-bi-invariant continuous
positive definite normalized functions on the subgroups ${\Cal
G}(n)$. \qed
\endproclaim

Note, however, that Theorem 3.2, due to its second claim, provides us
with more detailed information on the approximation process than
Theorem 3.5.

\remark{Remark\/ \rm3.6} Note that Claim 2 in the proof of Theorem 3.2 can
also be deduced from Doob's theorem on convergence of (reversed)
martingales (see \cite{D, Chapter VII, Theorem 4.2}). Indeed, denote
by ${\Cal B}_n$ the $\sigma$-algebra of all ${\Cal K}(n)$-invariant
Borel subsets of the space ${\Cal X}$. We have ${\Cal
B}_1\supseteq{\Cal B}_2\supseteq\dots$, so that, by Doob's theorem,
as $n\to\infty$, the conditional expectation $E(\psi\mi{\Cal B}_n)$
of the bounded Borel function $\psi$ converges almost everywhere,
with respect to $\nu$, to a function, which is the conditional
expectation $E(\psi\mi{\Cal B}_{\infty})$, where ${\Cal
B}_{\infty}=\bigcap{\Cal B}_n$. To derive Claim 2, we only need to
show that $E(\psi\mi{\Cal B}_n)=\psi_n$ almost everywhere with
respect to $\nu$, $n=1,2,\dots$. By definition of conditional
expectation, this means that $\psi_n$ is ${\Cal B}_n$-measurable and
the following condition holds:
$$
A\in{\Cal B}_n\implies\int_A\psi(x\)\nu(dx)
=\int_A\psi_n(x\)\nu(dx)\.\tag3.19
$$

Since $\psi_n$ is a ${\Cal K}(n)$-invariant Borel function, it is
${\Cal B}_n$-measurable. Further, since the action of ${\Cal K}$ on
${\Cal X}$ is jointly continuous, the function
$(u,x)\mapsto\psi(u\cdot x)$ is continuous on ${\Cal K}\times{\Cal
X}$, hence is a Borel function on ${\Cal K}(n)\times A$. By Fubini's
theorem,
$$
\int_{{\Cal K}(n)}\bigg(\int_A\psi(u\cdot x\)\nu(dx)\bigg)m_n(du)=
\int_A\bigg(\int_{{\Cal K}(n)}\psi(u\cdot x\)m_n(du)\bigg)\nu(dx)\.\tag3.20
$$
The left-hand side of (3.20) is equal to $\int_A\psi(x\)\nu(dx)$,
because $A$ and $\nu$ are ${\Cal K}(n)$-invariant, whereas the
right-hand side is equal to $\int_A\psi_n(x\)\nu(dx)$ by the
definition of $\psi_n$.
\endremark

\head\S4.\enspace Main Theorem \endhead

We shall deal with a sequence $\{M_n\in{\Cal M}_n\}$, $n=1,2,\dots$,
of orbital measures. By $f_n$ we denote the characteristic function
of $M_n$; recall that
$$
f_n(A)=\int_{\Omega_n}e^{i\op{tr}(AB)}M_n(dB),\qquad A\in H(n),\tag4.1
$$
where $\Omega_n$ stands for the $U(n)$-orbit that carries $M_n$. Let
$\Lambda(n)=(\lambda_1(n),\dots,\lambda_n(n))$ be the common
spectrum of all the matrices $B\in\Omega_n$. Then $\Omega_n$ may be
specified as the orbit containing the diagonal matrix
$\op{diag}\Lambda(n)$ with diagonal entries
$(\lambda_1(n),\dots,\lambda_n(n))$. The eigenvalues
$\lambda_1(n),\dots,\lambda_n(n)$ may be arranged in any order; it
will be convenient for us to separate the positive and the negative
eigenvalues and to regard $\Lambda(n)$ as a double sequence formed by
positive and negative eigenvalues, respectively, written in
decreasing order of their absolute values:
$$
\Lambda(n)=(\Lambda'(n),\Lambda''(n)),\tag4.2
$$
where
$$
\aligned
\Lambda'(n)&=(\lambda'_1(n)\ge\lambda'_2(n)\ge\cdots\ge0),\\
\Lambda''(n)&=(\lambda''_1(n)\le\lambda''_2(n)\le\cdots\le0)\.
\endaligned\tag4.3
$$
The possible zero values may be included either in $\Lambda'(n)$ or
in $\Lambda''(n)$; in fact, we prefer to view both $\Lambda'(n)$ and
$\Lambda''(n)$ as infinite sequences with a finite number of nonzero
terms.

\proclaim{Theorem 4.1 \rm(Main Theorem)} Let $\{M_n\in{\Cal
M}_n\}$ be an infinite sequence of orbital measures defined by a
sequence $\{\Omega_n\subset H(n)\}$ of $U(n)$-orbits. For
$n=1,2,\dots$, pick a matrix $B_n$ from the orbit $\Omega_n$ and
write the collection $\Lambda(n)$ of its eigenvalues as a double
sequence\/ \rom{(4.2). (}Note that $\Lambda(n)$ does not depend on
the choice of $B_n$.\rom)

\rom{(i)} Suppose that the following limits exist\/\rom:
$$
\gather
x'_k=\lim_{n\to\infty}\frac{\lambda'_k(n)}{n}\ge0,\quad
x''_k=\lim_{n\to\infty}\frac{\lambda''_k(n)}{n}\le0,\qquad
k=1,2,\dots,\tag4.4\\
\gamma_1=\lim_{n\to\infty}\frac1n
\sum_k(\lambda'_k(n)+\lambda''_k(n))=
\lim_{n\to\infty}\frac1n\op{tr}B_n,\tag4.5\\
\tilde\gamma_2=\lim_{n\to\infty}\frac1{n^2}
\sum_k((\lambda'_k(n))^2+(\lambda''_k(n))^2)=
\lim_{n\to\infty}\frac1{n^2}\op{tr}(B_n^2)\.
\tag4.6
\endgather
$$
Then the measures $M_n$ weakly converge to an ergodic measure
$M\in{\Cal M}$ with the multiplicative characteristic function $f$
defined by
$$
f(A)=\prod_{a\in\op{Spec}(A)}F(a),\qquad A\in H(\infty),\tag4.7
$$
where $\op{Spec}\(A)$ stands for the collection
$\{a_1,a_2,\dots,0,0,\dots\}$ of the eigenvalues of $A$ and
$$
F(a)=e^{i\gamma_1a-\gamma_2a^2\!/2}
\prod_k\frac{e^{-ix'_ka}}{1-ix'_ka}
\prod_k\frac{e^{-ix''_ka}}{1-ix''_ka},\qquad a\in{\Bbb R};\tag4.8
$$
here $\gamma_1$ is given by {\rm(4.5)}, the parameters $x'_k\ge0$ and
$x''_k\le0$ are given by {\rm(4.4)}, and, finally,
$$
\gamma_2=\tilde\gamma_2-\sum_k((x'_k)^2+(x''_k)^2),\tag4.9
$$
where $\tilde\gamma_2$ is given by\/ \rom{(4.6)}.

\rom{(ii)} Conversely, if the measures $M_n$ weakly converge to a
probability measure on $H$, then the limits {\rm (4.4)--(4.6)} do exist.
\endproclaim

\remark{Comment} It follows from the definition of the parameters
$x'_k$, $x''_k$, and $\tilde\gamma_2$ that
$\sum((x'_k)^2+(x''_k)^2)\le\tilde\gamma_2$, so that $\gamma_2\ge0$.
\endremark

\proclaim{Corollary 4.2} Let $\{M_n\in{\Cal M}_n\}$ be a sequence of
orbital measures that weakly converges to a Borel probability measure
$M$ on $H$. Then $M\in{\Cal M}$. \qed
\endproclaim

\demo{Derivation of Theorem\/ \rom{2.9} from Theorem\/ \rom{4.1}}
Let $F\in{\Cal F}_1$ and
let $M\in{\Cal M}$ be the corresponding ergodic measure. We
must show that $F$ is of the form (2.9). By Theorem 3.3, there exists
a sequence $\{M_n\in{\Cal M}_n\}$ that weakly converges to $M$. Next,
by claim (ii) of Theorem 4.1, the limits (4.4)--(4.6) exist.
Finally, by claim (i) of Theorem 4.1, $F$ is of the form (4.8) that
coincides with (2.9) up to a reordering of the points in
$x=(x_1,x_2,\dots)$ only; we recall (see Comment 1 to Proposition
2.8) that we may order the parameters $x_1,x_2,\dots$ in (2.9) in
any way. \qed
\enddemo

\demo{Outline of the proof of Theorem\/ \rm4.1} In \S5, we establish a
preliminary result -- we expand the orbital integral (4.1) into a
series of Schur polynomials. The proof of the theorem is given in
\S6; it is divided into three steps.

In Step 1, we check that under assumptions (4.4)--(4.6),
$$
f_n(\op{diag}\(a,0,0,\dots))\to F_{\gamma_1,\gamma_2,x}(a),\qquad
a\in{\Bbb R},\tag4.10
$$
where $\op{diag}\(\cdots)$ stands for a diagonal matrix and
$F_{\gamma_1,\gamma_2,x}$ is given by (4.8) or, that is the same, by
(2.9).

In Step 2, we generalize this to arbitrary diagonal matrices:
$$
f_n(\op{diag}\(a_1,\dots,a_k,0,0,\dots))\to
\prod_{i=1}^k F_{\gamma_1,\gamma_2,x}(a_i),\qquad
(a_1,\dots,a_k)\in{\Bbb R}^k.\tag4.11
$$
This proves claim (i).

Finally, in Step 3, using a simple trick, we show that conditions
(4.4)--(4.6) are indeed necessary (claim (ii)).
\enddemo

Let us emphasize that to prove Theorem 2.9 only, one could avoid Step~2.
However, we need this step to characterize the convergent
sequences $\{M_n\}$.

\remark{Remark\/ \rm4.3} To check (4.10) or, more generally, (4.11), we
consider the Taylor series decomposition at zero for the left-hand
side and show that its coefficients tend, as $n\to\infty$, to the
corresponding Taylor coefficients for the right-hand side. Note that
these coefficients are nothing but the moments of the measures.

Moreover, it follows from claim (ii) and the proof of (i) that
whenever a sequence $\{M_n\in{\Cal M}_n\}$ weakly converges to a
probability measure $M$, the moments of $M_n$ must tend to the
moments of $M$. Thus, in our situation, weak convergence $M_n\to M$
is always controlled by moments---a fact that is not at all
evident a priori.

Such a ``moment method'' also works in allied classification
problems, related to characters of $U(\infty)$ (see Vershik--Kerov
\cite{VK2}) and spherical functions of $GL(\infty,{\Bbb C})$ (see
Nessonov \cite{N}). We conjecture that it can be used for all
families of classical symmetric spaces. (About spherical functions on
infinite-dimensional symmetric spaces, see Olshanski\u\i\ \cite{O5} and
Pickrell \cite{Pi1}.)
\endremark

\remark{Remark\/ \rm4.4} Assume that in the spectrum
$\Lambda(n)=(\lambda_1(n),\dots,\lambda_n(n))$ there are at most $k$
nonzero eigenvalues, where $k$ does not depend on $n$. Then we can
prove claim (i) of Theorem 4.1 directly, i.e., without using
moments. Indeed, let
$$
\lim_{n\to\infty}\frac{\lambda_i(n)}n=x_i,\quad1\le i\le n,\qquad
\lambda_i(n)=0,\quad i>k\.\tag4.12
$$
Then the measure $M_n$ can be viewed as the image of the normalized
Haar measure of the group $U(n)$ under the mapping
$$
U(n)\ni u\mapsto u^*\op{diag}\(\lambda_1(n),\dots,\lambda_k(n),0,
\dots,0\)u=B\in H(n)\subset H\.\tag4.13
$$
The matrix $B$ can be rewritten as follows:
$$
B=(\Xi(n))^*X(n\)\Xi(n),\tag4.14
$$
where $\Xi(n)$ denotes the $k\times n$ matrix formed by the first $k$
rows of the matrix $u\in U(n)$ multiplied by the scalar $\sqrt{n}$,
and
$$
X(n)=\op{diag}\bigg(\frac{\lambda_1(n)}n,\dots,\frac{\lambda_k(n)}n
\bigg).\tag4.15
$$
Now let $n\to\infty$. Then $X(n)\to X:=\op{diag}\(x_1,\dots,x_k)$,
because of (4.12). Further, fix $m=1,2,\dots$ and regard the
$k\times m$ matrix formed by the first $m$ columns of $\Xi(n)$ as a
random matrix variable (with respect to the Haar measure of $U(n)$).
It is well-known that the limit distribution of this matrix is given
by the Gaussian product measure $\omega^{k\times m}$ (product of $km$
copies of $\omega$, cf.~(2.12)), where $\omega$ stands for the
Gaussian measure on ${\Bbb C}$ specified in Corollary 2.7. This fact
is proved, e.g., in Olshanski \cite{O5, Lemma 5.3}. It follows that
the measures $M_n$ weakly converge, as $n\to\infty$, to the measure
$M_{0;x_1,\dots,x_k}\in{\Cal M}$ defined in Remark 2.10.
\endremark

\remark{Remark\/ \rm4.5} In general, for a convergent sequence of orbital
measures, the eigenvalues in the spectrum $\Lambda(n)$ must grow
linearly in $n$ as $n\to\infty$. But if the limiting ergodic measure
is Gaussian (that is, $x\equiv0$), then the order of growth of
the eigenvalues becomes equal to $\sqrt{n}$. Example: let
$$
\Lambda(n)=(\underbrace{\sqrt{\gamma n},\dots,\sqrt{\gamma n}}_{[n/2]},\;
\underbrace{-\sqrt{\gamma n},\dots,-\sqrt{\gamma n}}_{[(n+1)/2]}\,),
$$
then the limiting measure is the Gaussian measure $M_{0,\gamma,0}$.
\endremark

\head\S5.\enspace Expanding spherical functions into series of Schur
polynomials\endhead

In this section, we fix $n=1,2,\dots$. Let $(a_1,\dots,a_n)\in{\Bbb
C}^n$, $(\lambda_1,\dots,\lambda_n)\in{\Bbb R}^n$, and suppose
$$
A=\op{diag}\(a_1,\dots,a_n),\qquad\Lambda=\op{diag}\(\lambda_1,\dots,
\lambda_n)\tag5.1
$$
are the corresponding diagonal matrices. We fix $\Lambda$ and deal
with the $U(n)$-orbit $\Omega\subset H(n)$ passing through $\Lambda$.
Let $f_{\Lambda}$ stand for the characteristic function of the
invariant probability measure supported by $\Omega$. Then
$$
f_{\Lambda}(A)=\int_{U(n)}e^{i\op{tr}(Au\Lambda u^{-1})}\,du,\tag5.2
$$
where $du$ is the normalized Haar measure on $U(n)$.

Since $f_{\Lambda}$ is an entire function of
$(a_1,\dots,a_n)\in{\Bbb C}^n$, it admits an everywhere convergent
Taylor series expansion. But since $f_{\Lambda}(a_1,\dots,a_n)$ is
also symmetric in $a_1,\dots,a_n$, it is more convenient to rewrite
this Taylor series as a series of Schur polynomials $s_{\mu}$,
$$
f_{\Lambda}(A)=f_{\Lambda}(a_1,\dots,a_n)=
\sum_{\mu}c_{\mu}s_{\mu}(a_1,\dots,a_n),\tag5.3
$$
where $\mu$ ranges over the set of all Young diagrams with at most
$n$ rows. (About Schur polynomials, see, e.g., Macdonald \cite{M}.)

We shall use some standard notation concerning Young diagrams:
$\mu\vdash m$ means that $\mu$ is a partition of $m$ (i.e., $m$
equals $|\mu|$, the number of boxes in $\mu$), $\ell(\mu)$ is the
number of (nonzero) rows in $\mu$, $(p,q)\in\mu$ denotes the box of
$\mu$ lying on the intersection of $p$\<th row and $q$\<th column,
$\dim\mu$ is the dimension of the irreducible representation of the
symmetric group $S(m)$, $m=|\mu|$, that corresponds to the diagram
$\mu$.

\proclaim{Theorem 5.1} The coefficients in\/ \rom{(5.3)} are given by
the following formula\/\rom:
$$
c_{\mu}=\prod_{(p,q)\in\mu}\frac1{n+q-p}\cdot s_{\mu}
(i\lambda_1,\dots,i\lambda_n)=i^{|\mu|}\prod_{(p,q)\in\mu}\frac
1{n+q-p}\cdot s_{\mu}(\lambda_1,\dots,\lambda_n)\.\tag5.4
$$
\endproclaim
\demo{Proof} {\it Step\/} 1. Let $\pi$ be an irreducible
representation of $U(n)$, $\op{Dim}\pi$ be its dimension, and $\chi$
be its normalized character:
$$
\chi(g)=\frac{\op{tr}\pi(g)}{\op{Dim}\pi},\qquad g\in U(n)\.
$$
It is well known that $\chi$ satisfies the following functional
equation:
$$
\int_{U(n)}\chi(guhu^{-1})\,du=\chi(g\)\chi(h),\qquad g,h\in U(n)\.\tag5.5
$$

Now let us take as $\pi$ the irreducible polynomial representation
with highest weight $\mu=(\mu_1,\dots,\mu_n)$. Then
$$
\chi(g)=\frac{s_{\mu}(\op{Spec}\(g))}{\op{Dim}_n\mu},\qquad g\in
U(n),\tag5.6
$$
where, given a $n\times n$ matrix $g$, $\op{Spec}\(g)$ stands for the
collection $(z_1,\dots,z_n)$ of its eigenvalues, and $\op{Dim}_n\mu$
denotes the dimension of the representation $\pi$ of the group $U(n)$.

Substituting (5.6) into (5.5), we obtain
$$
\int_{U(n)}s_{\mu}(\op{Spec}\(guhu^{-1}))\,du=
\frac1{\op{Dim}_n\mu}\,s_{\mu}(\op{Spec}\(g)\)s_{\mu}(\op{Spec}\(h))\.\tag5.7
$$
By analytic continuation, this formula holds for arbitrary $n\times
n$ complex matrices $g$ and $h$.

{\it Step\/} 2. Let us come back to the orbital integral (5.2), where
$A$ and $\Lambda$ are given by (5.1). We may write
$$
e^{i\op{tr}(Au\Lambda u^{-1})}=\sum_{m\ge0}
\frac1{m!}\,p_1^m(\op{Spec}\(iAu\Lambda u^{-1})),\tag5.8
$$
where $p_1(z_1,\dots,z_n)=z_1+\cdots+z_n$ is the first power sum.

Recall the well-known identity (see, e.g., Macdonald \cite{M, chapter
I, (7.8)}):
$$
p_1^m(z_1,\dots,z_n)=\sum\Sb\mu\vdash m\\
\ell(\mu)\le n\endSb \dim\mu\cdot s_{\mu}(z_1,\dots,z_n)\.\tag5.9
$$
Using it, we may rewrite (5.8) as follows:
$$
e^{i\op{tr}(Au\Lambda u^{-1})}=\sum_{m\ge0}\sum\Sb\mu\vdash m\\
\ell(\mu)\le n\endSb
\frac{\dim\mu}{m!}\,s_{\mu}(\op{Spec}\(iAu\Lambda u^{-1}))\.\tag5.10
$$
Integrating both sides of (5.10) over $u\in U(n)$ and applying
the functional equation (5.7) with $g=A$, $h=i\Lambda$, we obtain
$$
\int_{U(n)}e^{i\op{tr}(Au\Lambda u^{-1})}\,du
=\sum_{m\ge0}\sum\Sb\mu\vdash m\\
\ell(\mu)\le n\endSb \frac{\dim\mu}{m!\op{Dim}_n\mu}\,
s_{\mu}(a_1,\dots,a_n\)s_{\mu}(i\lambda_1,\dots,i\lambda_n)\.\tag5.11
$$

{\it Step\/} 3. For a box $(p,q)\in\mu$, let $h(p,q)$ denote the
corresponding hook length. Recall the well-known formulas
$$
\dim\mu=m!\prod_{(p,q)\in\mu}\frac1{h(p,q)},
\qquad\op{Dim}_n\mu=\prod_{(p,q)\in\mu}\frac{n+q-p}{h(p,q)},
$$
see, e.g., \cite{M, Chapter I, \S3, Example 4}. It follows
$$
\frac{\dim\mu}{m!\op{Dim}_n\mu}=\prod_{(p,q)\in\mu}\frac1{n+q-p}\.
$$
Substituting this into (5.11), we obtain (5.4). \qed
\enddemo

(After work on this paper was completed, we learned
that the argument presented above was used much earlier by James, see his
survey paper \cite{J, (60)}.)

\proclaim{Corollary 5.2} The orbital integral\/ \rom{(5.2)} is given by
the following explicit formula\/\rom:
$$
f_{\Lambda}(A)=\int_{U(n)}e^{i\op{tr}(Au\Lambda u^{-1})}\,du=
\frac{(n-1)!\cdots0!\det\>[e^{ia_j\lambda_k}]_{j,k=1}^n}
{V(a_1,\dots,a_n\)V(i\lambda_1,\dots,i\lambda_n)}\.\tag5.12
$$
Here $V(\cdots)$ is the Vandermonde determinant,
$$
V(z_1,\dots,z_n)=\prod_{1\le j<k\le n}(z_j-z_k),
$$
and the coordinates
$a_1,\dots,a_n$, as well as $\lambda_1,\dots,\lambda_n$, are
assumed to be pairwise distinct.
\endproclaim

\demo{Proof} The statement of Theorem 5.1 may be rewritten as follows:
$$
f_{\Lambda}(A)=(n-1)!\cdots0!\sum\Sb\mu\\
\ell(\mu)\le n\endSb
\frac{s_{\mu}(a_1,\dots,a_n\)s_{\mu}(i\lambda_1,\dots,i\lambda_n)}
{(\mu_1+n-1)!\,(\mu_2+n-2)!\cdots\mu_n!}\.
$$
Using the determinant formula for Schur polynomials, we obtain
$$
\split
&\frac{V(a_1,\dots,a_n\)V(i\lambda_1,\dots,i\lambda_n)}{(n-1)!\cdots0!}\,
f_{\Lambda}(A)\\
&\qquad=\sum\Sb\mu\\
\ell(\mu)\le n\endSb \frac{\det\>[a_j^{\mu_k+n-k}]_{j,k=1}^n\cdot
\det\>[(i\lambda_j)^{\mu_k+n-k}]_{j,k=1}^n}
{(\mu_1+n-1)!\,(\mu_2+n-2)!\cdots\mu_n!}\\
&\qquad=\sum_{m_1>\cdots>m_n\ge0} \frac{\det\>[a_j^{m_k}]_{j,k=1}^n\cdot
\det\>[(i\lambda_j)^{m_k}]_{j,k=1}^n} {m_1!\,m_2!\cdots m_n!}\.
\endsplit
$$

Since the numerator in the latter expression is symmetric in
$m_1,\dots,m_n$ and vanishes if some of these numbers are equal, we
may drop the assumption $m_1>\cdots>m_n$. Then we obtain
$$
\split
&\frac{V(a_1,\dots,a_n\)V(i\lambda_1,\dots,i\lambda_n)}
{(n-1)!\cdots0!}\,f_{\Lambda}(A)\\
&\qquad=\frac1{n!}\sum_{m_1,\dots,m_n\ge0}
\frac{\det\>[a_j^{m_k}]_{j,k=1}^n\cdot \det\>[(i\lambda_j)^{m_k}]_{j,k=1}^n}
{m_1!\,m_2!\cdots m_n!}\\
&\qquad=\frac1{n!}\sum_{m_1,\dots,m_n\ge0}\frac1{m_1!\cdots m_n!}
\sum_{\sigma,\tau}\op{sgn}\(\sigma)\op{sgn}(\tau)
(a_{\sigma_1}\cdot i\lambda_{\tau_1})^{m_1}\cdots
(a_{\sigma_n}\cdot i\lambda_{\tau_n})^{m_n} ,
\endsplit
$$
where $\sigma=(\sigma_1,\dots,\sigma_n)$ and
$\tau=(\tau_1,\dots,\tau_n)$ range over all the permutations of
$1,\dots,n$.

Changing the order of the summation, we obtain
$$
\frac1{n!}\sum_{\sigma,\tau}\op{sgn}\(\sigma)\op{sgn}\(\tau\)
e^{ia_{\sigma_1}\lambda_{\tau_1}}\cdots e^{ia_{\sigma_n}\lambda_{\tau_n}}
=\det\>[e^{ia_j\lambda_k}]_{j,k=1}^n,
$$
which completes the proof. \qed
\enddemo

\remark{Remark\/ \rm5.3} Formula (5.12) is well known and can be
proved in a number of different ways. For instance, it can be
obtained from the Gelfand--Naimark calculation of spherical functions
on $GL(n,{\Bbb C})$ (see Gelfand--Naimark \cite{GN}) by a
passage to the limit.
Another way consists in using the radial parts of invariant
differential operators. Yet another method can be found in \cite{BGV,
Section 7.5} (we are grateful to Michel Duflo for the latter
reference). Note also that writing the above calculations in reverse
order, we can deduce Theorem 5.1 from formula (5.12).
\endremark

\proclaim{Corollary 5.4} Let us substitute
$A=\op{diag}\(a,0,\dots,0)$ into\/ \rom{(5.2)}. Then
$$
f_{\Lambda}(\op{diag}\(a,0,\dots,0))=\int_{U(n)}e^{ia(u\Lambda
u^{-1})_{11}}\,du=\sum_{m\ge0}\frac{h_m(i\lambda_1,\dots,i\lambda_n)}
{n(n+1)\cdots(n+m-1)}\,a^m,\tag5.13
$$
where $h_m$ is the $m$\<th complete symmetric function.
\endproclaim

\demo{Proof} Indeed, note that $s_{\mu}(a,0,\dots,0)$ vanishes
unless $\mu=(m)$, where $m=0,1,\dots$. Thus the summation is really
taken over the diagrams $\mu=(m)$ only. For these diagrams, $s_{\mu}$
reduces to $h_m$ and the product $\prod_{(p,q)\in\mu}(n+q-p)$ turns
into $n(n+1)\cdots(n+m-1)$, so that we obtain (5.13). \qed
\enddemo

\remark{Remark\/ \rm5.5} The expansion (5.3) of theorem 5.1 may be
written in the following form, which emphasizes the symmetry between
$A$ and $\Lambda$,
$$
f_{\Lambda}(A)=\sum\Sb\mu\\
\ell(\mu)\le n\endSb i^{|\mu|}\bigg(\prod_{(p,q)\in\mu}\frac1{n+q-p}\bigg)
s_{\mu}(\lambda_1,\dots,\lambda_n\)s_{\mu}(a_1,\dots,a_n)\.\tag5.14
$$
The symmetry can already be seen from (5.2).
\endremark

\head\S6.\enspace Proof of theorem 4.1 \endhead

We start with the following simple claim.\footnote{See also Remark
6.2 at the end of this section.}

\proclaim{Proposition 6.1} Let $f,f_1,f_2,\dots$ be analytic
functions on ${\Bbb R}^k$ satisfying the following conditions\/\rom:
\roster
\item"(i)" $f_1,f_2,\dots$ are positive definite and normalized at
the origin\/\rom;
\item"(ii)" the Taylor coefficients of $f_n$ at the origin tend, as
$n\to\infty$, to the corresponding Taylor coefficients of $f$\rom;
\item"(iii)" The Taylor decomposition of $f_n$ converges absolutely
and uniformly on $n$ in a neighborhood of the origin that does not
depend on $n$.
\endroster

Then $f_n\to f$ uniformly on compact subsets of ${\Bbb R}^k$.
\endproclaim

\demo{Proof} This is a standard exercise. It is clear that $f_n\to f$
in a neighborhood of the origin. Consider the probability measures
$\nu_1,\nu_2,\dots$ on ${\Bbb R}^k$ that correspond to the functions
$f_1,f_2,\dots$ by Bochner's theorem. By Paul L\'evy's classical
continuity theorem (see, e.g., Shiryaev's textbook \cite{Shi,
chapter III, \S3}), the measures $\nu_n$ weakly converge to a
probability measure $\nu$. Let $\tilde f$ stand for the
characteristic function of $\nu$. Then $f_n\to\tilde f$ uniformly on
compact sets of ${\Bbb R}^k$. Further, (ii) and (iii) imply that
$\tilde f$ is analytic in a neighborhood of the origin, and this
implies that $\tilde f$ is analytic on the whole space ${\Bbb R}^k$.
Therefore, $\tilde f=f$ and our claim follows. \qed
\enddemo

We proceed to prove the theorem. The notation of \S4 is 
maintained.

{\it Step\/} 1. Let us fix a sequence $\{\Lambda(n)\}$ such that the
limits (4.4), (4.5), and (4.6) exist. Let us abbreviate
$$
f_n(a)=f_n(\op{diag}\(a,0,0,\dots)),\qquad a\in{\Bbb R}\.\tag6.1
$$
Note that $f_n(a)$ is the characteristic function of the probability
measure $\nu_n:=M_n^{(1)}$, the image of the measure $M_n$ under the
projection $\theta_1\:H(n)\to H(1)={\Bbb R}$.

The purpose of this step is to prove that
$$
\lim_{n\to\infty}f_n(a)=F_{\gamma_1,\gamma_2,x}(a),\qquad a\in{\Bbb R},
\tag6.2
$$
uniformly on bounded sets in ${\Bbb R}$. (Recall that $\gamma_1$ is
given by (4.5), $\gamma_2$ is given by (4.6) and (4.9), and
$x=(x'_1,x'_2,\dots;x'_1,x''_2,\dots)$, where the latter
parameters are defined by (4.4). The function
$F_{\gamma_1,\gamma_2,x}$ is given by (4.8) or, that is the same, by
(2.9).)

To do this, let us expand both the sides of (6.2) into Taylor series:
$$
\align
f_n(a)&=\sum_{m\ge0}c_m^{(n)}a^m,\tag6.3\\
F_{\gamma_1,\gamma_2,x}(a)&=\sum_{m\ge0}c_m^{(\infty)}a^m.\tag6.4
\endalign
$$
By Proposition 6.1, it is enough to verify the following two claims:
First,
$$
\lim_{n\to\infty}c_m^{(n)}=c_m^{(\infty)},\qquad m=0,1,2,\dots\.\tag6.5
$$
Second, the series (6.3) converges absolutely and uniformly on $n$
in a sufficiently small neighborhood of the origin, i.e.,
$$
|c_m^{(n)}|\le C_1C_2^m,\qquad m=0,1,2,\dots,\tag6.6
$$
where the constants $C_1>0$, $C_2>0$ do not depend on $n$.

By (5.13), we have
$$
\split
c_m^{(n)}&=\frac{h_m(i\lambda_1(n),\dots,i\lambda_n(n))}
{n(n+1)\cdots(n+m-1)}\\
&=\frac{n^m}{n(n+1)\cdots(n+m-1)}\,h_m\bigg(\frac{i\lambda_1(n)}n,\dots,
\frac{i\lambda_n(n)}n\bigg).
\endsplit\tag6.7
$$
It is clear that in both claims, (6.5) and (6.6), we may replace
$c_m^{(n)}$ by
$$
\tilde c_m^{(n)}:=h_m\bigg(\frac{i\lambda_1(n)}n,\dots,
\frac{i\lambda_n(n)}n\bigg). \tag6.8
$$

Further, instead of the series $\sum\tilde c_m^{(n)}a^m$ and $\sum
c_m^{(\infty)}a^m$, it is more convenient to deal with their
logarithms $\ln\(\sum\tilde c_m^{(n)}a^m)$ and $\ln\(\sum
c_m^{(\infty)}a^m)$, respectively.

Now recall a well-known identity from the theory of symmetric
functions,
$$
\sum_{m\ge0}h_m(\,\bdot\,\)a^m=\exp\bigg(\sum_{m\ge1}
p_m(\,\bdot\,)\,\frac{a^m}m\bigg),\tag6.9
$$
where $p_m(\,\bdot\,)$ are the power sum symmetric functions (see, e.g.,
\cite{M}). It follows from (6.8) and (6.9) that
$$
\ln\bigg(\sum_{m\ge0}\tilde c_m^{(n)}a^m\bigg)=\sum_{m\ge1}
p_m\bigg(\frac{i\lambda_1(n)}n,\dots,\frac{i\lambda_n(n)}n\bigg)
\frac{a^m}m\.\tag6.10
$$
On the other hand, by definition (2.9) of the function
$F_{\gamma_1,\gamma_2,x}$,
$$
\ln\bigg(\sum_{m\ge0}c_m^{(\infty)}a^m\bigg)=i\gamma_1a-\frac12\,\gamma_2
a^2+\sum_{m\ge2}p_m(ix)\,\frac{a^m}m,\tag6.11
$$
where
$$
p_m(ix)=i^m p_m(x)=i^m\sum_{k=1}^{\infty}((x'_k)^m+(x''_k)^m)\.\tag6.12
$$
Note that the sum in the right-hand side of (6.12) is convergent,
because, due to assumption (4.6), we have $p_2(x)<\infty$ (see the
Comment to Theorem 4.1).

Thus, our first claim reduces to the existence of the limits
$$
\align
\lim_{n\to\infty}p_1\bigg(\frac{i\lambda_1(n)}n,\dots,
\frac{i\lambda_n(n)}n\bigg)&=i\gamma_1,\tag6.13\\
\lim_{n\to\infty}p_2\bigg(\frac{i\lambda_1(n)}n,\dots,
\frac{i\lambda_n(n)}n\bigg)&=-\gamma_2+p_2(ix),\tag6.14\\
\lim_{n\to\infty}p_m\bigg(\frac{i\lambda_1(n)}n,\dots,
\frac{i\lambda_n(n)}n\bigg)&=p_m(ix),\qquad m\ge3,\tag6.15
\endalign
$$
and our second claim reduces to an estimate of the form
$$
\bigg|p_m\bigg(\frac{i\lambda_1(n)}n,\dots,\frac{i\lambda_n(n)}n\bigg)
\bigg|\le C'_1(C'_2)^m,\qquad m\ge3,\tag6.16
$$
where $C'_1>0$, $C'_2>0$ are some constants not depending on $n$.

Clearly, (6.13) is just the assumption (4.5). Further, (6.14)
immediately follows from (4.6) and (4.9). Indeed, by (4.6),
$$
\lim_{n\to\infty}p_2\bigg(\frac{i\lambda_1(n)}n,\dots,
\frac{i\lambda_n(n)}n\bigg)=-\lim_{n\to\infty}p_2
\bigg(\frac{\lambda_1(n)}n,\dots,\frac{\lambda_n(n)}n\bigg)
=-\tilde\gamma_2\.
$$
Now, by (4.9),
$$
-\tilde\gamma_2=-\gamma_2-\sum_k((x'_k)^2+(x''_k)^2)=-\gamma_2+p_2(ix)\.
$$
Let us verify (6.15). Using the notation (4.2), (4.3), we have
$$
\align
p_m\bigg(\frac{i\lambda_1(n)}n,\dots,\frac{i\lambda_n(n)}n\bigg)
&=\sum_{r\ge1}\bigg(\frac{i\lambda'_r(n)}n\bigg)^m+
\sum_{r\ge1}\bigg(\frac{i\lambda''_r(n)}n\bigg)^m,\tag6.17\\
p_m(ix)&=\sum_{r\ge1}(ix'_r)^m+\sum_{r\ge1}(ix''_r)^m.\tag6.18
\endalign
$$
To deduce (6.15) from (4.4), it suffices to show that both sums
in the right-hand side of (6.17) converge absolutely and uniformly
on $n$. Let us examine the first sum (for the second one the reasoning
is just the same). Since $m\ge3$ and
$\lambda'_1(n)\ge\lambda'_2(n)\ge\cdots$, we have for $N=1,2,\dots$
$$
\split
\bigg|\sum_{r\ge N}\bigg(\frac{i\lambda'_r(n)}n\bigg)^m\bigg|
&=\sum_{r\ge N}\bigg(\frac{\lambda'_r(n)}n\bigg)^m
\le\bigg(\frac{\lambda'_N(n)}n\bigg)^{m-2}
\sum_{r\ge N}\bigg(\frac{\lambda'_r(n)}n\bigg)^2\\
&\le\frac{\lambda'_N(n)}n p_2
\bigg(\frac{\lambda_1(n)}n,\dots,\frac{\lambda_n(n)}n\bigg).
\endsplit\tag6.19
$$
By (4.6),
$p_2(\lambda_1(n)/n,\dots,\lambda_n(n)/n)$
remains bounded as $n\to\infty$. Finally, since $x'_N\to0$ as
$N\to\infty$ and since $\lambda'_N(n)/n\to x'_N$ as
$n\to\infty$, the value of $\lambda'_N(n)/n$ may be made
arbitrarily small provided first $N$ and then $n$ are chosen large
enough.

Let us verify (6.16). For $m\ge3$
$$
\split
&\bigg|p_m\bigg(\frac{i\lambda_1(n)}n,\dots,\frac{i\lambda_n(n)}n\bigg)
\bigg|\\
&\qquad\le p_2\bigg(\frac{\lambda_1(n)}n,\dots,\frac{\lambda_n(n)}n\bigg)
\cdot\sup_{1\le r\le n}\bigg(\frac{\lambda_r(n)}n\bigg)^{m-2}\\
&\qquad\le p_2\bigg(\frac{\lambda_1(n)}n,\dots,\frac{\lambda_n(n)}n\bigg)
p_2\bigg(\frac{\lambda_1(n)}n,\dots,\frac{\lambda_n(n)}
n\bigg)^{(m-2)/2}.
\endsplit\tag6.20
$$
Since $p_2(\lambda_1(n)/n,\dots,\lambda_n(n)/n)$
remains bounded as $n\to\infty$, we obtain (6.16).

This completes Step 1.

{\it Step\/} 2. The purpose of this step is to prove that under the
same assumptions as in Step 1, we have a more general result: for any
fixed $k=1,2,\dots$,
$$
\lim_{n\to\infty}f_n(\op{diag}\(a_1,\dots,a_k,0,\dots,0))=
\prod_{p=1}^k F_{\gamma_1,\gamma_2,x}(a_p)\tag6.21
$$
uniformly on bounded subsets in ${\Bbb R}^k$.

Let us abbreviate
$$
f_n(a_1,\dots,a_k)=f_n(\op{diag}\(a_1,\dots,a_k,0,\dots,0))\tag6.22
$$
and note that $f_n(a_1,\dots,a_k)$ is again the characteristic
function of some probability measure $\nu_n$ on ${\Bbb R}^k$. Namely,
$\nu_n$ is the radial part of the measure $M_n^{(k)}=\theta_k(M_n)$
on $H(k)$ with respect to the projection
$$
H(k)\ni A\mapsto\op{Spec}\(A)\in{\Bbb R}^k.\tag6.23
$$

Our reasoning will be similar to that of Step 1. We expand both
sides of (6.21) into multidimensional Taylor series. However, since
these are symmetric functions of $(a_1,\dots,a_k)$, we prefer to
rewrite the Taylor series as series of Schur polynomials
$s_{\mu}(a_1,\dots,a_k)$, where $\mu$ ranges over the set of all
Young diagrams with $\ell(\,\bdot\,)\le k$:
$$
\align
f_n(a_1,\dots,a_k)&=\sum c_{\mu}^{(n)}s_{\mu}(a_1,\dots,a_k), \tag6.24\\
\prod_{p=1}^kF_{\gamma_1,\gamma_2,x}(a_p)&=\sum c_{\mu}^{(\infty)}s_{\mu}
(a_1,\dots,a_k)\. \tag6.25
\endalign
$$

As in Step 1, applying Proposition 6.1, we reduce our problem to
verifying the following two claims:

First,
$$
\lim_{n\to\infty}c_{\mu}^{(n)}=c_{\mu}^{(\infty)}\quad
\text{for any $\mu$ with $\ell(\mu)\le k$},\tag6.26
$$
and, second,
$$
|c_{\mu}^{(n)}|\le C_1C_2^{|\mu|},\tag6.27
$$
where $C_1>0$, $C_2>0$ are some constants not depending on $n$.

By (5.4),
$$
\split
c_{\mu}^{(n)}&=\prod_{(p,q)\in\mu}\frac1{n+q-p}
\cdot s_{\mu}(i\lambda_1(n),\dots,i\lambda_n(n))\\
&=\prod_{(p,q)\in\mu}\frac n{n+q-p}\cdot s_{\mu}
\bigg(\frac{i\lambda_1(n)}n,\dots,\frac{i\lambda_n(n)}n\bigg).
\endsplit\tag6.28
$$
Thus, in both claims, (6.26) and (6.27), we may replace
$c_{\mu}^{(n)}$ by
$$
\tilde c_{\mu}^{(n)}:=s_{\mu}
\bigg(\frac{i\lambda_1(n)}n,\dots,\frac{i\lambda_n(n)}n\bigg).\tag6.29
$$

Let us prove (6.26) with $c_{\mu}^{(n)}$ replaced by $\tilde
c_{\mu}^{(n)}$.

Recall the Jacobi--Trudi identity expressing the Schur functions in
terms of the complete symmetric functions (see, e.g., \cite{M,
Chapter I, (3.4)}):
$$
s_{\mu}=\det\>[h_{\mu_i-i+j}]_{i,j=1}^k,\qquad\ell(\mu)\le k\.\tag6.30
$$
It follows from (6.8), (6.29), and (6.30) that
$$
\tilde c_{\mu}^{(n)}=\det\>[\tilde c_{\mu_i-i+j}^{(n)}]_{i,j=1}^k\.\tag6.31
$$
By Step 1, $\tilde c_m^{(n)}\to c_m^{(\infty)}$ as $n\to\infty$, so
that
$$
\lim_{n\to\infty}\tilde
c_{\mu}^{(n)}=\det\>[c_{\mu_i-i+j}^{(\infty)}]_{i,j=1}^k\.\tag6.32
$$
Further, it is well known that for an arbitrary formal series
$\sum_{m\ge0}c_ma^m$ with $c_0=1$, we have
$$
\prod_{p=1}^k\bigg(\sum_{m=0}^{\infty}c_ma_p^m\bigg)=\sum\Sb\mu\\
\ell(\mu)\le k\endSb
\det\>[c_{\mu_i-i+j}]_{i,j=1}^k s_{\mu}(a_1,\dots,a_k)\.\tag6.33
$$
Applying (6.33) to the left-hand side of (6.25), we conclude that
$$
c_{\mu}^{(\infty)}=\det\>[c_{\mu_i-
i+j}^{(\infty)}]_{i,j=1}^k=\lim_{n\to\infty} \tilde c_{\mu}^{(n)}.\tag6.34
$$

Thus, we have verified the first claim. As for the second one, an
estimate of type (6.27) for $\tilde c_{\mu}^{(n)}$ follows at once
from (6.31) and from the estimate (6.6) for $\tilde c_m^{(n)}$, proved in
Step 1.

This completes Step 2.

It follows that $f_n\to f$ (where $f$ is given by (4.7)) uniformly
on compact subsets of $H(\infty)$. Clearly, $f$ is a continuous
positive definite normalized function on $H(\infty)$. By Proposition
1.1, it is the characteristic function of a $U(\infty)$-invariant
Borel probability measure $M$ on the space $H$, invariant under the
action of $U(\infty)$. Then, by Proposition 3.4, $M_n$ weakly
converges to $M$. Since $f$ is multiplicative, $M$ is ergodic. Thus,
we have verified claim (i) of Theorem 4.1.

{\it Step\/} 3. Let us fix a sequence $\{\Lambda(n)\}$, where
$\Lambda(n)=(\lambda_1(n),\dots,\lambda_n(n))$, and let $f_n(a)$,
$a\in{\Bbb R}$, be defined as in Step 1. Let us assume that
$$
\lim_{n\to\infty}f_n(a)=f(a),\qquad a\in{\Bbb R},\tag6.35
$$
uniformly on bounded subsets in ${\Bbb R}$, where $f(a)$ is a
function on ${\Bbb R}$ (of course, $f$ is automatically continuous).
We shall prove that then $\{\Lambda(n)\}$ must satisfy the
assumptions (4.4)--(4.6) of Theorem 4.1(i).

Suppose first that
$$
\sup_n\bigg\{p_2\bigg(\frac{\lambda_1(n)}n,\dots,\frac{\lambda_n(n)}n\bigg)
+\bigg(p_1\bigg(\frac{\lambda_1(n)}n,\dots,\frac{\lambda_n(n)}n\bigg)
\bigg)^2\bigg\} <\infty\.\tag6.36
$$
Then, given an infinite subset $N\subseteq\{1,2,\dots\}$, there
exists a possibly smaller infinite subset $N'\subseteq N$ such that
the limits (4.4), (4.5), and (4.6) exist provided $n$ goes to
infinity inside $N'$. Then, by Step 1, $f_n\to
F_{\gamma_1,\gamma_2,x}$, as $n\to\infty$ inside $N'$, so that
$F_{\gamma_1,\gamma_2,x}=f$. For any other $N$ and $N'$, the
parameters $\gamma_1$, $\gamma_2$, and $x$ will be the same, because
they are uniquely determined by the function itself, see Comment 4
after Proposition 2.8. It follows that the limits (4.4)--(4.6) exist
as $n$ ranges over the set of all natural numbers.

Suppose now that (6.36) does not hold. We shall show that this leads
to a contradiction with the initial assumption (6.35).

Indeed, since the expression $\{\cdots\}$ in (6.36) is a homogeneous
function of $\Lambda(n)$, we can choose an infinite subset
$N\subseteq\{1,2,\dots\}$ and a sequence of positive numbers
$\{\varepsilon_n\mid n\in N\}$ such that $\lim_{n\in
N}\varepsilon_n=0$ and
$$
\lim_{n\in N}\bigg\{p_2\bigg(\varepsilon_n\frac{\lambda_1(n)}n,\dots,
\varepsilon_n \frac{\lambda_n(n)}n\bigg)+
\bigg(p_1\bigg(\varepsilon_n\frac{\lambda_1(n)}n,\dots,
\varepsilon_n\frac{\lambda_n(n)}n\bigg)\bigg)^2\bigg\}=1\.\tag6.37
$$
Then, replacing $N$ by a smaller infinite subset $N'$, we can arrange
so that for the sequence $\{\varepsilon_n\Lambda(n)\}$, the limits
(4.4)--(4.6) will exist provided $n$ goes to infinity inside $N'$.
Moreover, at least one of the corresponding parameters $\gamma_1$,
$\tilde\gamma_2$ will be nonzero.

Note that the effect of multiplying $\Lambda(n)$ by $\varepsilon_n$
is the same as that of multiplying $a$ by $\varepsilon_n$. Thus, by
Step 1,
$$
\lim_{n\in N'}f_n(\varepsilon_n a)=F_{\gamma_1,\gamma_2,x},\qquad a
\in{\Bbb R}\.\tag6.38
$$
Since at least one of the parameters $\gamma_1$, $\tilde\gamma_2$ of
the function $F_{\gamma_1,\gamma_2,x}$ is nonzero, it follows from
the definition of this function that it is not equal identically to
1. But since $F_{\gamma_1,\gamma_2,x}$ is analytic, the same is true
in an arbitrarily small neighborhood of the point $a=0$. Then, comparing
(6.38) with (6.35), we arrive at a contradiction.

Thus, we have verified claim (ii) of Theorem 4.1. \qed

\remark{Remark\/ \rm6.2} Note that the estimates (6.6) and (6.16) are not
necessary to assert the convergence of the functions $f_n$. Indeed,
Proposition 6.1 may be replaced by the following stronger claim:

{\it Let $f_1,f_2,\dots$ and $f$ be smooth positive definite
functions on ${\Bbb R}^k$, normalized at the origin. Expand them into
Taylor series at the origin and assume that each Taylor coefficient
of $f_n$ tends, as $n\to\infty$, to the corresponding coefficient of
$f$. Finally, assume that the moment problem defined by the
coefficients of $f$ has a unique solution\/ \rom(the latter condition is
satisfied, e.g., if $f$ is analytic\/\rom).

Then the sequence $(f_n)$ converges to $f$ uniformly on compact
subsets of ${\Bbb R}^n$.}
\endremark

By virtue of this claim, the verification of the uniform convergence of
the Taylor expansions may be omitted.

\head\S7.\enspace Total positivity\endhead

\definition{Definition 7.1} Let $\varphi(t)$ be a real nonnegative
measurable function on ${\Bbb R}$. Then $\varphi$ is said to be a
{\it totally positive\/} function if for $n=1,2,\dots$
$$
\det\>[\varphi(t_i-s_j)]_{i,j=1}^n\ge0\quad\text{for any
$t_1<\dots<t_n$ and $s_1<\dots<s_n$}\.\tag7.1
$$

It will be convenient for us to include in the definition the
following additional assumption: $\varphi$ is summable and
$\int\varphi(t)\,dt=1$, i.e., $\varphi(t)\,dt$ is a probability measure
on ${\Bbb R}$. (In \cite{S2}, functions satisfying both
conditions are called {\it P\'olya frequency functions}, the
second condition is in fact not restrictive, see \cite{S2, lemma 4}.)
\enddefinition

\proclaim{Proposition 7.2 \rm(Schoenberg \cite{S2, p.~341,
Lemma~5})} The set of totally positive functions is stable under
convolution.
\endproclaim

\demo{Proof} For two summable functions $\varphi$ and $\psi$, the
convolution $\varphi\ast\psi$ is correctly defined and the following
formula is readily verified:
$$
\multline
\det\>[(\varphi\ast\psi)(t_i-s_j)]_{i,j=1}^n\\
=\frac1{n!}\int_{{\Bbb R}^n}\det\>[\varphi(t_i-u_k)]_{i,k=1}^n\cdot
\det\>[\psi(u_k-s_j)]_{k,j=1}^n\,du_1\cdots du_n\.
\endmultline\tag7.2
$$

Now suppose $\varphi$ and $\psi$ are totally positive and let
$t_1<\cdots<t_n$, $s_1<\cdots<s_n$. Then the integrand in (7.2) is
nonnegative for all (pairwise distinct) $u_1,\dots,u_n$, because
both determinants have the same sign, equal to that of
$\prod_{k>l}(u_k-u_l)$.\qed
\enddemo

\proclaim{Proposition 7.3 \rm(Schoenberg \cite{S2, p.~335 and
p.~343})} The densities of the normal and exponential distributions,
$$
\psi_{\gamma}(t)=\frac1{\sqrt{2\pi\gamma}}\,e^{-t^2\!/(2\gamma)},
\qquad\gamma>0,\tag7.3
$$
and
$$
\varphi_y(t)=\cases y^{-1}e^{-y^{-1}t},&t\ge0,\\
0,&t<0,\endcases\qquad y>0,\tag7.4
$$
are totally positive.\qed
\endproclaim

Note that the result remains true after the shift $t\mapsto
t+\op{const}$ of the argument or the change of sign $t\mapsto-t$.

Note also that $\{\psi_{\gamma}\}$ and $\{\varphi_y\}$ are
one-parametric semigroups with respect to the convolution product.

\proclaim{Theorem 7.4 \rm(Schoenberg's theorem on totally positive
functions, see Schoenberg \cite{S2}, Karlin \cite{K})} The Fourier
transforms $\widehat\varphi$ of totally positive functions $\varphi$
are just the functions $F_{\gamma_1,\gamma_2,x}$, defined in {\rm
(2.9)}, where at least one of the parameters $\gamma_2,x_1,x_2,\dots$ is
nonzero.\qed
\endproclaim

This fundamental result shows that a totally positive function
$\varphi$ is a convolution product
$\varphi_0\ast\varphi_1\ast\varphi_2\ast\cdots$, where
$\varphi_0(t)\,dt$ is a normal distribution and $\varphi_k(t)\,dt$,
$k=1,2,\dots$, are, up to transformations $t\mapsto\pm t+\op{const}$,
exponential distributions.

(We have to exclude $\gamma_2=x_1=x_2=\cdots=0$, because the inverse
Fourier transform of the corresponding function $F$ is a Dirac
measure.)

Comparing Theorem 2.9 and Theorem 7.4, we obtain the following
correspondence $M\leftrightarrow\varphi$ between the ergodic measures
$M\in{\Cal M}$ (except the Dirac measures on scalar matrices) and the
totally positive functions $\varphi$ on the real line:
$$
\theta_1(M)(dt)=\varphi(t)\,dt
$$
(recall that the mapping $\theta_1$ assigns to a matrix $B\in H$ its
matrix element $B_{11}\in{\Bbb R}$). In other words, this means that
the distribution of the random variable $B_{11}$ with respect to the
probability distribution $M$ on the matrices $B\in H$ is given by the
density $\varphi$.

The easy part of Theorem 7.4 consists in verifying the fact that
$F_{\gamma_1,\gamma_2,x}=\widehat\varphi$ with a totally positive
$\varphi$. This is done by making use of Proposition 7.2 and
Proposition 7.3 and an evident passage to the
limit (cf.~Proposition 2.8).

The hard part of Theorem 7.4 is to prove that the Fourier transform
$\widehat\varphi$ of any totally positive function $\varphi$ is of
the form (2.9). Our purpose is to show that this claim is equivalent
to Theorem 2.9.

\definition{Definition 7.5 \rm(Karlin \cite{K, p.~49})} A real smooth
nonnegative function $\varphi(t)$ on ${\Bbb R}$,
$\int\varphi(t)\,dt=1$, is called {\it extended totally positive\/} if
$$
\det\>[\varphi^{(i-1)}(v_j)]_{i,j=1}^n\ge0,\qquad
n=1,2,\dots,\;v_1>\dots>v_n\.\tag7.5
$$
\enddefinition

\proclaim{Proposition 7.6} \rom{(i)} Any smooth totally positive
function $\varphi$ is extended totally positive.

\rom{(ii)} Conversely, if $\varphi$ is an extended totally positive
function, then the function $\varphi\ast\psi_{\gamma}$, where
$\psi_{\gamma}$ was defined by {\rm (7.3)}, is totally positive for
any $\gamma>0$.
\endproclaim

\demo{Proof} (i) By definition of total positivity, for any pairwise
distinct real $t_1,\dots,t_n$ and any $s_1<\cdots<s_n$,
$$
\prod_{p>q}(t_p-t_q)^{-1}\cdot\det\>[\varphi(t_i-s_j)]_{i,j=1}^n\ge 0\.
\tag7.6
$$
Putting $s_1=-v_1,\dots,s_n=-v_n$ and letting $t_1,\dots,t_n\to0$
in (7.6), we obtain (7.5).

(ii) By Theorem 2.1 in Karlin \cite{K, p.~50}, if $\varphi$ verifies
strict inequalities in (7.5), then it is totally positive. So it suffices
to prove that if we replace $\varphi$ by
$\varphi\ast\psi_{\gamma}$, then the inequalities in (7.5) become
strict.

For $n=1$ this is evident, because $\varphi$ is nonnegative and not
identically equal to zero whereas $\psi_{\gamma}$ is strictly
positive. For $n>1$ this argument is generalized as follows.

First, remark that the functions $\varphi,\varphi',\varphi'',\dots$
are linearly independent. Indeed, if this is not true, then $\varphi$
satisfies a linear differential equation with constant coefficients,
whence $|\varphi(t)|\to\infty$ as $t\to\infty$ or $t\to-\infty$. But
this contradicts the assumption $\varphi\in L^1({\Bbb R})$.

Next, substitute $\psi=\psi_{\gamma}$ into formula (7.2) and repeat
the argument used in the proof of (i). Then we obtain, for any
$s_1<\cdots<s_n$,
$$
\multline
\det\>[(\varphi\ast\psi_{\gamma})^{(i-1)}(-s_j)]_{i,j=1}^n\\
=\frac1{n!}\int_{{\Bbb R}^n}\det\>[\varphi^{(i-1)}(-u_k)]_{i,k=1}^n
\det\>[\psi_{\gamma}(u_k-s_j)]_{k,j=1}^n\,du_1\cdots du_n\.
\endmultline\tag7.7
$$
Arguing just as in the proof of Proposition 7.2, we see that the
integrand in (7.7) is nonnegative.

Finally, as noted in Schoenberg \cite{S2, p.~336}, the second
determinant in the integrand is nonzero provided $u_1,\dots,u_n$ are
pairwise distinct. On the other hand, since $\varphi,\varphi',\dots$
are linearly independent, the first determinant in the integrand does
not vanish for certain $(u_1,\dots,u_n)\in{\Bbb R}^n$, hence on an
open subset of ${\Bbb R}^n$. We conclude that the integrand is
everywhere nonnegative and strictly positive on an open subset, so
that the integral (7.7) is strictly positive. \qed
\enddemo

Note that Proposition 7.6 corresponds to a part of Pickrell's proof
in \cite{Pi2, p.~154--155}. There it is claimed that an analytic
extended totally positive function is totally positive; however, the
arguments are too sketchy and seem to be incomplete. To avoid this
difficulty, we modified the claim somewhat and used a trick suggested
by Boyer's paper \cite{Bo, p.~218}.

Note also that not all totally positive functions are smooth: for
instance, the function (7.4) is not smooth at 0. Thus, the class of
extended totally positive functions, as defined above, does not
coincide with the class of totally positive functions (although
the two classes are very close, as is seen from Proposition 7.6). For
this reason, to use property (7.5), we must first smooth totally
positive functions.

The next theorem is Pickrell's main calculation in \cite{Pi2,
pp.~154--155}. It is simple but instructive. For completeness and
for reader's convenience, we give the proof (which is presented here in
slightly more detail than in \cite{Pi2}).

\proclaim{Theorem 7.7 \rm(Pickrell \cite{Pi2, pp.~154--155})} Let
$\varphi$ be a smooth nonnegative function on ${\Bbb R}$,
$\int\varphi(t)\,dt=1$, and let $F=\widehat\varphi$ be its Fourier
transform. Then $F$ belongs to the class ${\Cal F}_1$ \rom(see
definition\/ \rom{(2.2))} if and only if $\varphi$ is extended totally
positive.
\endproclaim
\demo{Proof} Given $n=1,2,\dots$ and $A\in H(n)$, denote by
$a_1,\dots,a_n$ the eigenvalues of $A$ and put
$$
f_n(A)=F(a_1)\cdots F(a_n);
$$
this is a continuous $U(n)$-invariant function on $H(n)$. Let us fix
$n$ and show that positive definiteness of $f_n$ is equivalent to
condition (7.5).

By Bochner's theorem, $f_n$ is positive definite if and only if its
inverse Fourier transform is a measure. This condition is equivalent
to the following one: for any function $\Psi\ge0$ from the Schwartz
space $\Cal S(H(n))$,
$$
\langle f_n,\widehat\Psi\rangle:=\int_{H(n)}f_n(A\)\ov{\widehat\Psi(A)}\,
dA\ge0\.\tag7.8
$$
Since $f_n$ is $U(n)$-invariant, one may assume $\Psi$ is
$U(n)$-invariant too.

Let $D(n)$ denote the subspace of diagonal matrices in $H(n)$. We
identify $D(n)$ with ${\Bbb R}^n$ and write elements of $D(n)$ as
$\op{diag}\(a_1,\dots,a_n)$ where $(a_1,\dots,a_n)\in{\Bbb R}^n$. It
is well known and easily verified that the radial part of the
Lebesgue measure on $H(n)$ with respect to the action of $U(n)$ is
the measure
$$
\op{const}V^2(a_1,\dots,a_n)\,da_1\cdots da_n,\qquad\op{const}>0\tag7.9
$$
on $D(n)$, where
$$
V(a_1,\dots,a_n)=\prod_{p<q}(a_p-a_q)\.
$$
It follows that
$$
\multline
\langle f_n,\widehat\Psi\rangle=\op{const}\int_{D(n)}
V(a_1,\dots,a_n\)\widehat\varphi(a_1)\cdots\widehat\varphi(a_n)\\
\times\ov{V(a_1,\dots,a_n\)\widehat\Psi(\op{diag}\(a_1,\dots,a_n))}\,
da_1\cdots da_n\.
\endmultline\tag7.10
$$
Put
$$
\theta(t_1,\dots,t_n)=V(t_1,\dots,t_n\)\Psi(\op{diag}\(t_1,\dots,t_n))\.
\tag7.11
$$
We shall show that the integral (7.10) is equal, up to a positive
factor, to
$$
\int_{{\Bbb R}^n}\det\>[\varphi^{(j-1)}(t_k)]_{j,k=1}^n\theta
(t_1,\dots,t_n)\,dt_1\cdots dt_n\.\tag7.12
$$
Indeed, the function
$V(a_1,\dots,a_n\)\widehat\varphi(a_1)\cdots\widehat\varphi(a_n)$ is
the Fourier transform of
$$
V\bigg(i\frac{\partial}{\partial t_1},\dots,
i\frac{\partial}{\partial t_n}\bigg)\cdot\varphi(t_1)\cdots\varphi(t_n)
=i^{n(n-1)/2}\det\>[\varphi^{(n-j)}(t_k)]_{j,k=1}^n\.\tag7.13
$$
On the other hand, using the $U(n)$-invariance of $\Psi$, we have
$$
\split
&V(a_1,\dots,a_n\)\widehat\Psi(\op{diag}\(a_1,\dots,a_n))\\
&\qquad=V(a_1,\dots,a_n)\int_{T\in
H(n)}e^{i\op{tr}(\op{diag}\(a_1,\dots,a_n)T)}\Psi(T)\,dT\\
&\qquad=V(a_1,\dots,a_n)\int_{H(n)}
\bigg(\int_{U(n)}e^{i\op{tr}(\op{diag}\(a_1,\dots,a_n)uTu^{-1})}du\bigg)
\Psi(T)\,dT.
\endsplit\tag7.14
$$
Using formula (5.12) for the interior integral and again applying
formula (7.8) for the radial part of the Lebesgue measure, we see
that (7.14) is equal, up to a positive factor, to
$$
(-i)^{n(n-1)/2}\int_{\Bbb R^n}\det\>[e^{ia_jt_k}]_{j,k=1}^n\cdot
\theta(t_1,\dots,t_n)\,dt_1\cdots dt_n,\tag7.15
$$
where $\theta$ was defined in (7.11). Developing the determinant and
using the fact that $\theta(t_1,\dots,t_n)$ is antisymmetric with
respect to permutations of $t_1,\dots,t_n$, we conclude that (7.15)
is equal to
$$
n!\(-i)^{n(n-1)/2}\hat\theta(a_1,\dots,a_n)\.\tag7.16
$$

Now (7.13) and (7.16) imply that the integral (7.10) is equal, up to
a positive factor, to
$$
\split
&(-1)^{n(n-1)/2}\int_{{\Bbb R}^n}\det\>[\varphi^{(n-j)}(t_k)]_{j,k=1}^n\cdot
\theta(t_1,\dots,t_n)\,dt_1\cdots dt_n\\
&\qquad=\int_{{\Bbb R}^n}\det\>[\varphi^{(j-1)}(t_k)]_{j,k=1}^n\cdot
\theta(t_1,\dots,t_n)\,dt_1\dots dt_n\.
\endsplit\tag7.17
$$
Thus, we have verified (7.12).

Since the integrand in (7.12) is symmetric with respect to
permutations of $t_1,\dots,t_n$, we see that (7.12) is equal, up to
a positive factor, to
$$
\int_{t_1>\cdots>t_n}\det\>[\varphi^{(j-1)}(t_k)]_{j,k=1}^n\cdot
\theta(t_1,\dots,t_n)\,dt_1\dots dt_n\.\tag7.18
$$
Recall that $\theta$ is given by (7.11), where $\Psi$ is an
$U(n)$-invariant nonnegative function from the Schwartz space, and
remark that $V(t_1,\dots,t_n)>0$ in the domain $t_1>\cdots>t_n$. It
follows that (7.18) is nonnegative for any such $\Psi$ if and only if
$\varphi$ is extended totally positive. \qed
\enddemo

\proclaim{Corollary 7.8} Theorem\/ \rom{2.9 (}classification of ergodic
measures $M\in{\Cal M}$ or of functions $F\in{\Cal F}_1$\rom) and
Schoenberg's Classification Theorem\/ \rom{7.4} can be derived one from
another.
\endproclaim

\demo{Proof} Let us show that Theorem 2.9 implies Theorem 7.4. Let
$\varphi$ be a totally positive function. We must prove that the
Fourier transform $F=\widehat\varphi$ is of the form (2.9), where at
least one of the parameters $\gamma_2,x_1,x_2,\dots$ is nonzero.

First suppose $\varphi$ is smooth. Then, by Proposition 7.6 (i),
$\varphi$ is extended totally positive. Next, by Theorem 7.7,
$F\in{\Cal F}_1$, and, finally, by Theorem 2.9, $F$ is of the form
(2.9).

The general case can be reduced to that of a smooth $\varphi$ as
follows. We again use the Gaussian totally positive function
$\psi_{\gamma}$ (see (7.3)) to smooth $\varphi$. Then, for any
$\gamma>0$, we have a smooth (even analytic) totally positive
function $\varphi\ast\psi_{\gamma}$. This implies that
$F(a\)e^{-\gamma a^2\!/2}$ is of the form (2.9) for any
$\gamma>0$, whence $F$ itself is of this form.

It should be added that $F$ cannot be equal to $F_{\gamma_1,0,0}$,
because we know that the inverse Fourier transform of $F$ is a
function and not a Dirac measure. So at least one of the parameters
$\gamma_2,x_1,x_2,\dots$ is nonzero.

The inverse implication is verified similarly, by making use of
Proposition 7.6 (ii). \qed
\enddemo

Thus, our proof of Theorem 2.9 leads to a new proof of Schoenberg's
Theorem~7.4.

\head\S8.\enspace Totally positive functions as limits of splines\endhead

After reading the preliminary version \cite{OV} of the present paper,
Andre\u\i\ Okounkov remarked that the one-dimensional projections of
orbital measures coincide with the so-called fundamental splines (=
$B$-splines) whose limits were studied in an important paper by Curry
and Schoenberg \cite{CS}. The purpose of this section is to briefly
discuss the relationship between the results of \cite{CS} and our results.

We start by stating some classical facts used in \cite{CS}.

Fix real numbers $t_1<\dots<t_n$, called the {\it knots\/}
($n\ge3$). There exists a (unique) function
$M_{n-1}(t)=M_{n-1}(t;t_1,\dots,t_n)$ on ${\Bbb R}$ such that:
\roster
\item"(i)" On each open interval determined by adjacent knots, $M_{n-1}(t)$
is a polynomial of degree $n-2$.
\item"(ii)" $M_{n-1}(t)$ vanishes when $t<t_1$ or $t>t_n$.
\item"(iii)" $M_{n-1}(t)$ has $n-3$ continuous derivatives at each knot.
\item"(iv)" $\int M_{n-1}(t)\,dt=1$.
\endroster

In \cite{CS}, the function $M_{n-1}(t)$ is called the {\it
fundamental spline\/} (with knots $t_1,\dots,t_n$). Another term,
used in the modern literature, is {\it $B$-spline.}

The fundamental spline $M_{n-1}(t)$ is given by the following
explicit formula:
$$
M_{n-1}(t;t_1,\dots,t_n)=(n-1)\sum_{k=1}^n\frac{(\max(t_k-t,0))^{n-2}}
{\prod_{i\ne k}(t_k-t_i)}\.\tag8.1
$$

Let $\sigma_{n-1}$ denote the standard $(n-1)$-dimensional simplex,
$$
\sigma_{n-1}=\{(p_1,\dots,p_n)\mid
0\le p_1,\dots,p_n\le1,\;p_1+\cdots+p_n=1\}\subset{\Bbb R}^n,\tag8.2
$$
and let $\xi$ denote the affine functional on the simplex taking
values $t_1,\dots,t_n$ at its vertices. Then $M_{n-1}(t)$ coincides
with the density of the image under $\xi$ of the Lebesgue measure on
$\sigma_{n-1}$, so normalized that the volume of $\sigma_{n-1}$ is
equal to 1. This implies, in particular, that $M_{n-1}(t)$ is
nonnegative, so that $M_{n-1}(t)\,dt$ is a probability measure on
${\Bbb R}$.

Using a passage to the limit, one easily extends the definition of
$M_{n-1}(t)$ to the case when some of the knots coincide.

For further properties of the functions $M_{n-1}(t)$ and for proofs
of the facts mentioned above, see \cite{CS} or, e.g., Babenko's
textbook \cite{Ba}.

Now we are in a position to state the main result of Curry and
Schoenberg \cite{CS, Theorem 6}:

\proclaim{Theorem 8.1} Consider the class of probability measures on
${\Bbb R}$ which can be obtained as weak limits of measures of the
form $M_{n-1}(t;t_1,\dots,t_n)\,dt$, where $n\to\infty$ and the knots
$t_1,\dots,t_n$ depend on $n$. Then the characteristic functions of
the measures of this class are exactly those given by formula\/
\rom{(4.8)}.
\endproclaim

We shall briefly describe the method of proof used in \cite{CS}.
Given $t_1,\dots,t_n$, put
$$
F_n(a)=\int_{-\infty}^\infty\bigg(1-\frac{iat}n\bigg)^{-n}
M_{n-1}(t;t_1,\dots,t_n)\,dt\.\tag8.3
$$
Since
$$
\Big(1-\frac{iat}n\Big)^{-n}=e^{iat}\Big(1+O\Big(\frac1n\Big)\Big),\tag8.4
$$
$F_n(a)$ may be viewed as the `approximate Fourier transform' of
$M_{n-1}(t)$. Its advantage with respect to the ordinary Fourier
image of $M_{n-1}(t)$ is that it is given by a very simple expression,
namely
$$
F_n(a)=\prod_{k=1}^n\Big(1-\frac{iat_k}n\Big)^{-1}.\tag8.5
$$

Now note that $F_n(ia)^{-1}$ is a polynomial in $a$ with only
real zeros, and use a well-known theorem, due to Laguerre and
P\'olya, which describes the class of entire functions that can be
approximated by polynomials with real zeros: up to change of a
variable $a\mapsto ia$, these are exactly the reciprocals to
functions of type (4.8), see Hirschman--Widder \cite{HW}.

Note that this result on entire functions also plays an important
role in Schoenberg's classification of totally positive functions
(Theorem 7.4 above).

Comparing Theorem 8.1 with Theorem 7.4, we conclude that a
probability measure on ${\Bbb R}$ (distinct from a Dirac measure) can
be approximated by a sequence of fundamental splines with growing
number $n$ of knots if and only if it is given by a totally positive
density.

The following fact, remarked by Andre\u\i{} Okounkov, is crucial for our
discussion.

\proclaim{Proposition 8.2 \rm(A.~Yu.~Okounkov)} Let
$\lambda_1\le\dots\le\lambda_n$ be real numbers and $n\ge3$. Consider
the $U(n)$-orbit in $H(n)$ passing through the diagonal matrix
$\Lambda=\op{diag}\(\lambda_1,\dots,\lambda_n)$ and denote by
$\mu(dB)$ the corresponding orbital measure.

Then the image of $\mu$ under the projection $H(n)\ni B\mapsto
B_{11}\in{\Bbb R}$ coincides with the fundamental spline
$M_{n-1}(t;\lambda_1,\dots,\lambda_n)\,dt$.
\endproclaim

\demo{Proof} For $u\in U(n)$, we have
$$
(u\Lambda u^{-1})_{11}=\sum_{k=1}^n u_{1k}\lambda_k(u^{-1})_{k1}=
\sum_{k=1}^n|u_{1k}|^2\lambda_k\.\tag8.6
$$
When the matrix $u$ ranges over $U(n)$, its first row
$$
z=(z_1,\dots,z_n):=(u_{11},\dots,u_{1n})\tag8.7
$$
ranges over the unit sphere $S^{2n-1}\subset{\Bbb C}^n$, and under
the mapping $u\mapsto z$, the orbital measure $\mu$ projects onto the
normalized invariant measure on $S^{2n-1}$, which may be written as
$$
\op{const}\,\frac{d\(\op{Re}z_1\)d\(\op{Im}z_1)\cdots
d\(\op{Re}z_n\)d(\op{Im}z_n)}{d\(|z_1|^2+\cdots+|z_n|^2-1)}\.\tag8.8
$$

Further, under the mapping
$$
(z_1,\dots,z_n)\mapsto(|z_1|^2,\dots,|z_n|^2) =(p_1,\dots,p_n)\tag8.9
$$
of the sphere onto the simplex $\sigma_{n-1}$, the measure (8.8)
projects onto the measure
$$
\op{const}\,\frac{dp_1\cdots dp_n}{d\(p_1+\cdots+p_n-1)},\tag8.10
$$
which coincides with the normalized Lebesgue measure on the simplex.

Finally, the right-hand side of (8.6), which may be rewritten as
$\sum_{k=1}^np_k\lambda_k$, is just the value at the point
$(p_1,\dots,p_n)\in\sigma_{n-1}$ of the linear functional $\xi\:{\Bbb
R}^n\to{\Bbb R}$ taking values $\lambda_1,\dots,\lambda_n$ at the
vertices of the simplex. By a property of the spline function
$M_{n-1}(t)$ mentioned above, we conclude that the image of the
orbital measure $\mu$ under the mapping $B\mapsto B_{11}$ is equal to
$M_{n-1}(t;\lambda_1,\dots,\lambda_n)\,dt$. \qed
\enddemo

By virtue of Proposition 8.2, the one-dimensional projections of
orbital measures admit a very nice analytic interpretation in terms
of splines, and Curry--Schoenberg's result described in Theorem 8.1
turns out to be almost equivalent to the `one-dimensional part' of
our Theorem 4.1, that is, to the results of Steps 1 and 3 in \S6.
The difference in the statements is that Curry and Schoenberg do not
obtain necessary and sufficient conditions on the knots under which
a sequence of fundamental splines would be weakly convergent; they are only
interested in describing the limiting functions.

Theorems 7.4 and 8.1 together imply that totally positive functions
are exactly those functions which may be approximated by fundamental
splines with a growing number of knots. This fact seems to be highly
nontrivial, because the fundamental splines themselves are not totally
positive. We hope that the chain of relations traced in the present
paper furnishes a certain explanation of this phenomenon.

\remark{Remark\/ \rm8.3} Recall that the {\it Dirichlet distribution\/}
${\Cal D}(\theta_1,\dots,\theta_n)$ with parameters
$\theta_1>0,\dots,\theta_n>0$ is defined as the probability measure
on the $(n-1)$-dimensional simplex (8.2) whose density with respect
to the Lebesgue measure is given by
$$
\op{const}p_1^{\theta_1-1}\cdots p_n^{\theta_n-1},\tag8.11
$$
see Kingman \cite{Ki, Section 9.1}. Now in Proposition 8.2
let us replace the space $H(n)$ of $n\times n$ Hermitian matrices
by the space
of $n\times n$ real symmetric (respectively, quaternion Hermitian)
matrices. Then the one-dimensional projections of orbital measures
coincide with various one-dimensional projections of the Dirichlet
distribution ${\Cal D}(1/2,\dots,1/2)$ (respectively, of the
Dirichlet distribution ${\Cal D}(2,\dots,2)$).

More generally, one can consider one-dimensional projections of the
Dirichlet distribution ${\Cal D}(\theta,\dots,\theta)$ with
arbitrary parameter $\theta>0$. If $\theta=1,2,3,\dots$ then
one-dimensional projections of this distribution are the fundamental
splines with multiple knots: the multiplicity of each knot is equal
to $\theta$. For general $\theta>0$ there is no such
interpretation. However, for any $\theta$, using the `moment method',
one can still obtain an analog of Theorem 8.1. The limiting
measures will have a characteristic function of the following form
(cf.~(4.8)):
$$
F(a)=e^{i\gamma_1a-\gamma_2a^2\!/2}
\prod_k\frac{e^{-i\theta x'_ka}}{(1-ix'_ka)^\theta}
\prod_k\frac{e^{-i\theta x''_ka}}{(1-ix''_ka)^\theta},\qquad a\in{\Bbb R},
\tag8.12
$$
where the parameters are the same as in (4.8).

Finally, note that the parameter $\alpha=\theta^{-1}$ 
exactly corresponds to the
parameter that appears in the theory of Jack's symmetric functions (see
Stanley \cite{Sta} or the 2nd edition (1995) of Macdonald's book
\cite{M}). In particular, the $m$th moment of a one-dimensional
projection of the Dirichlet distribition 
${\Cal D}(\theta,\dots,\theta)$ is equal, up to a scalar factor, to
$P_m(t_1,\ldots,t_n;\theta^{-1})$, where $t_1,\ldots,t_n$ stand for the
parameters of the projection (i.e., the values of the corresponding affine
functional $\xi$ at the vertices of the simplex) and
$P_m(\;\cdot\;;\alpha)$, $m=1,2,\ldots$, are 
one-row Jack's symmetric functions with parameter $\alpha$.
\endremark

\enddocument

\Refs

\widestnumber\key{BGW}

\ref\key ASW
\by M. Aissen, I. J. Schoenberg,
and A. M. Whitney
\paper On the generating functions of
totally positive sequences, \rom I
\jour J. Analyse Math.
\vol 2
\yr 1952
\pages 93--103
\endref

\ref\key Ba
\by K. I. Babenko
\book Basic numerical analysis
\publ``Nauka''
\publaddr Moscow
\yr 1986
\lang Russian
\endref

\ref\key BCR
\by C. Berg, J. P. R. Christensen,
and P. Ressel
\book Harmonic analysis on semigroups.
Theory of positive definite and
related functions
\publ Springer-Verlag
\yr 1984
\endref

\ref\key BGV
\by N. Berline, E. Getzler,
and M. Vergne
\book Heat kernels and
Dirac operators
\bookinfo Grundlehren der math.
Wiss. 298
\publ Springer-Verlag
\yr 1992
\endref

\ref\key Bi
\by P. Billingsley
\paper Convergence of probability
measures
\publ Wiley
\publaddr New York etc.
\yr 1968
\endref

\ref\key Bo
\by R. P. Boyer
\paper Infinite traces of
AF-algebras and
characters of $U(\infty)$
\jour J. Operator Theory
\vol 9
\yr 1983
\pages 205--236
\endref

\ref\key CS
\by H. B. Curry and I. J. Schoenberg
\paper On P\'olya frequency
functions\/ \rom{IV}. The fundamental
spline functions and their limits
\jour J. Analyse Math.
\vol 17
\yr 1966
\pages 71--107
\endref

\ref\key D
\by J. L. Doob
\book Stochastic processes
\publ Wiley
\publaddr New York etc.
\yr 1953
\endref

\ref\key E1
\by A. Edrei
\paper On the generating functions of
totally positive
sequences\/ \rom{II}
\jour J. Analyse Math.
\vol 2
\yr 1952
\pages 104--109
\endref

\ref\key E2
\bysame
\paper On the generating function of a
doubly-infinite, totally positive sequence
\jour Trans. Amer. Math. Soc.
\vol 74
\yr 1953
\pages 367--383
\endref

\ref\key GN
\by I. M. Gelfand and M. A. Naimark
\book Unitary representations of
classical groups
\publaddr Moscow--Leningrad
\yr 1950
\lang Russian
\transl German transl.
\publ Akademie Verlag
\publaddr Berlin
\yr 1957
\endref

\ref\key HW
\by I. I. Hirschman and D. V. Widder
\book The convolution transform
\publaddr Princeton
\yr 1955
\endref

\ref\key I1
\by R. S. Ismagilov
\paper Linear representations of
groups of matrices with elements
from a normed field
\jour Izv. Akad. Nauk SSSR Ser. Mat.
\vol 33
\yr 1969
\pages 1296--1323
\transl English transl.
\jour Math. USSR-Izv.
\vol 3
\yr 1969
\pages 1219--1244
\endref

\ref\key I2
\bysame
\paper Spherical functions over
a normed field
whose residue field is infinite
\jour Funktsional. Anal. i Prilozhen.
\vol 4
\issue 1
\yr 1970
\pages 42--51
\transl English transl.
\jour Functional Anal. Appl.
\vol 4
\yr 1970
\issue 1
\pages 37--45
\endref

\ref\key J
\by A. T. James
\paper Distributions of matrix variates and latent roots derived from
normal samples
\jour Ann. Math. Stat.
\vol 35
\yr 1964
\pages 475--501
\endref

\ref\key K
\by S. Karlin
\book Total positivity
\vol I
\publ Stanford Univ. Press
\publaddr Stanford, CA
\yr 1968
\endref

\ref\key Ke
\by S. Kerov
\paper Gaussian limit for the
Plancherel measure of
the symmetric group
\jour C. R. Acad. Sci. Paris, S\'er. I
\vol 316
\yr 1993
\pages 303--308
\endref

\ref\key KV
\by S. V. Kerov and A. M. Vershik
\paper The characters of the
infinite symmetric group and
probability properties of the
Robinson--Schensted--Knuth algorithm
\jour SIAM J. Alg. Discr. Meth.
\vol 7
\yr 1986
\pages 116--124
\endref

\ref\key KOV
\by S. Kerov, G. Olshanski,
and A. Vershik
\paper Harmonic
analysis on the infinite symmetric group.
Deformation of the regular
representation
\jour C. R. Acad. Sci. Paris, S\'er. I
\vol 316
\yr 1993
\pages 773--778
\endref

\ref\key Ki
\by J. F. C. Kingman
\book Poisson processes
\bookinfo Oxford Studies in Probability
\vol 3
\publ Clarendon Press
\publaddr Oxford
\yr 1993
\endref

\ref\key M
\by I. G. Macdonald
\book Symmetric functions and Hall
polynomials
\publ Clarendon Press
\publaddr Oxford
\yr 1979
\endref

\ref\key Mu
\by R. J. Muirhead
\book Aspects of multivariate
statistical theory
\publ Wiley
\yr 1982
\endref

\ref\key N
\by N. I. Nessonov
\paper The complete classification of
representations of the group
$GL(\infty)$ containing the identity
representation of the unitary subgroup
\jour Mat. Sb.
\vol 130
\yr 1986
\pages 131--150
\transl English transl.
\jour Math. USSR-Sb.
\vol 5
\yr 1987
\pages 122--147
\endref

\ref\key Ok
\by A. Yu. Okounkov
\paper Thoma's theorem and
representations of
the infinite bisymmetric group
\jour Funktsional. Anal. i Prilozhen.
\vol 28
\issue 2
\yr 1994
\pages 31--40
\transl English transl.
\jour Functional Anal. Appl.
\vol 28
\yr 1994
\issue 2
\pages 100--107
\endref

\ref\key O1
\by G. I. Ol\mz shanski\u\i
\paper Unitary representations
of the infinite-dimensional
classical groups $U(p,\infty)$,
$SO(p,\infty)$, $Sp(p,\infty)$
and the corresponding motion groups
\jour Funktsional. Anal. i Prilozhen.
\vol 12
\issue 3
\yr 1978
\pages 32--44
\transl English transl.
\jour Functional Anal. Appl.
\vol 12
\yr 1978
\issue 3
\pages 185--195
\endref

\ref\key O2
\bysame
\paper Infinite-dimensional classical
groups of finite ${\Bbb R}$-rank:
description of representations and
asymptotic theory
\jour Funktsional. Anal. i Prilozhen.
\vol 18
\issue 1
\yr 1984
\pages 28--42
\transl English transl.
\jour Functional Anal. Appl.
\vol 18
\yr 1984
\issue 1
\pages 22--34
\endref

\ref\key O3
\bysame
\paper Unitary representations
of the group
$SO(\infty,\infty)$ as limits of
unitary representations of the
groups $SO(n,\infty)$ as $n\to\infty$
\jour Funktsional. Anal. i Prilozhen.
\vol 20
\issue 4
\yr 1986
\pages 46--57
\transl English transl.
\jour Functional Anal. Appl.
\vol 20
\yr 1986
\issue 4
\pages 292--301
\endref

\ref\key O4
\bysame
\paper Unitary representations of
$(G,K)$-pairs connected with
the infinite symmetric group
$S(\infty)$
\jour Algebra i Analiz
\vol 1
\issue 4
\yr 1989
\pages 178--209
\transl English transl.
\jour Leningrad Math. J.
\vol 1
\yr 1989
\pages 983--1014
\endref

\ref\key O5
\bysame
\paper Unitary representations of
infinite-dimensional pairs
$(G,K)$ and the formalism of R. Howe
\inbook  Representation of Lie
Groups and Related Topics
\eds A. M. Vershik and D. P. Zhelobenko
\bookinfo Advanced Studies in Contemporary
Mathematics
\vol 7
\publ Gordon and Breach Science Publishers
\publaddr New York etc.
\yr 1990
\pages 269--463
\endref

\ref\key O6
\bysame
\paper Irreducible unitary
representations
of the groups $U(p,q)$ sustaining
passage to the limit as
$q\to\infty$
\jour J. Soviet Math.
\vol 59
\issue 5
\yr 1992
\pages 1102--1107
\endref

\ref\key O7
\bysame
\paper Representations of
infinite-dimensional classical groups,
limits of enveloping algebras
and Yangians
\ed A.~A.~Kirillov
\bookinfo Topics in Representation
Theory. Advances in Soviet Math.
\vol 2
\publ Amer. Math. Soc.
\publaddr Providence, R.I.
\yr 1991
\pages 1--66
\endref

\ref\key OV
\by G. Olshanski and A. Vershik
\paper Ergodic unitarily invariant
measures on the space of infinite
Hermitian matrices
\jour Preprint UTMS 94--61 October 19,
University of Tokyo, Department
of Mathematical Sciences
\yr 1994
\endref

\ref\key Pa1
\by K. R. Parthasarathy
\book Probability measures on
metric spaces
\publ Academic Press
\publaddr New York--London
\yr 1967
\endref

\ref\key Pa2
\bysame
\book Introduction to probability
and measure
\yr 1980
\endref

\ref\key Ph
\by R. R. Phelps
\book Lectures on Choquet's theorem
\publ Van Nostrand
\yr 1966
\endref

\ref\key Pi1
\by D. Pickrell
\paper Separable representations
for automorphism
groups of infinite symmetric spaces
\jour J. Funct. Anal.
\vol 90
\yr 1990
\pages 1--26
\endref

\ref\key Pi2
\bysame
\paper Mackey analysis
of infinite classical
motion groups
\jour Pacific J. Math.
\vol 150
\yr 1991
\pages 139--166
\endref

\ref\key S1
\by I. J. Schoenberg
\paper Metric spaces and
completely monotone functions
\jour Ann. of Math.
\vol 39
\yr 1938
\pages 811--841
\endref

\ref\key S2
\by I. J. Schoenberg
\paper On P\'olya frequency
functions\/ \rom I. 
The totally positive functions and
their Laplace transforms
\jour Journal d'Analyse Math\'ematique
\vol 1
\yr 1951
\pages 331--374
\endref

\ref\key Sha
\by D. Shale
\paper Linear symmetries of
free boson fields
\jour Trans. Amer. Math. Soc.
\vol 103
\yr 1962
\pages 149--167
\endref

\ref\key Shi
\by A. N. Shiryaev
\book Probability
\publ``Nauka''
\publaddr Moscow
\yr 1980
\transl English transl.
\publ Springer-Verlag
\yr 1984
\endref

\ref\key Sta2
\by R. P. Stanley
\paper Some combinatorial properties of Jack symmetric functions
\jour Advances in Math.
\vol 77
\yr 1989
\pages 76--115
\endref

\ref\key SV
\by S. Stratila and D. Voiculescu
\paper A survey on
representations of unitary group
$U(\infty)$
\inbook Spectral Theory,
Banach Center Publications
\vol 8
\publ Polish Science Publ.
\publaddr Warsaw
\yr 1982
\pages 416--434
\endref

\ref\key T
\by E. Thoma
\paper Die unzerlegbaren,
positiv-definiten
Klassenfunktionen der abz\"ahlbar
unendlichen, symmetrischen Gruppe
\jour Math. Zeitschr.
\vol 85
\yr 1964
\pages 40--61
\endref

\ref\key V
\by A. M. Vershik
\paper Description of invariant
measures for the
actions of some infinite-dimensional
groups
\jour Dokl. Akad. Nauk SSSR
\vol 218
\yr 1974
\pages 749--752
\transl English transl.
\jour Soviet Math. Dokl.
\vol 15
\yr 1974
\pages 1396--1400
\endref

\ref\key VK1
\by A. M. Vershik and S. V. Kerov
\paper Asymptotic theory of
characters of the symmetric group
\jour Funktsional. Anal. i Prilozhen.
\vol 15
\issue 4
\yr 1981
\pages 17--27
\transl English transl.
\jour Functional Anal. Appl.
\vol 15
\yr 1981
\issue 4
\pages 246--255
\endref

\ref\key VK2
\bysame
\paper Characters and factor
representations of the infinite
unitary group
\jour Dokl. Akad. Nauk SSSR
\vol 267
\yr 1982
\pages 272--276
\transl English transl.
\jour Soviet Math. Dokl.
\vol 26
\yr 1982
\pages 570--574
\endref

\ref\key VK3
\bysame
\paper Asymptotics of the
largest and the typical dimensions of
the irreducible representations
of a symmetric group
\jour Funktsional. Anal. i Prilozhen.
\vol 19
\issue 1
\yr 1985
\pages 25--36
\transl English transl.
\jour Functional Anal. Appl.
\vol 19
\yr 1985
\pages 21--31
\issue 1
\endref

\ref\key VK4
\bysame
\paper The Grothendieck group of
the infinite symmetric group and
symmetric functions (with elements of
the theory of the $K_0$-functor of
$AF$-algebras)
\eds A. M. Vershik and D. P. Zhelobenko
\inbook Representation of Lie Groups
and Related
Topics. Advanced Studies in
Contemporary Mathematics
\vol 7
\publ Gordon and Breach Science Publishers
\publaddr New York etc.
\yr 1990
\pages 39--117
\endref

\ref\key Vo1
\by D. Voiculescu
\paper Repr\'esentations factorielles de type
\rom{II$_1$} de $U(\infty)$
\jour J. Math. Pures Appl.
\vol 55
\yr 1976
\pages 1--20
\endref

\ref\key Vo2
\bysame
\paper On extremal invariant functions of
positive type on certain groups
\jour INCREST Preprint Series Math.
\yr 1978
\endref

\endRefs
\bye